\documentclass{amsart}
\usepackage{amssymb}
\usepackage{amsmath,amscd}
\usepackage{color}
\usepackage{epsf,epsfig}
\usepackage{graphicx}
\theoremstyle{plain}

 \newtheorem{thm}{Theorem}[section]
 \newtheorem{cor}[thm]{Corollary}
 \newtheorem{lem}[thm]{Lemma}
 \newtheorem{prop}[thm]{Proposition}
 \theoremstyle{definition}
 \newtheorem{defn}{Definition}[section]
 \theoremstyle{remark}
 \newtheorem{rem}{Remark}[section]
 \newtheorem{example}{Example}[section]
 \numberwithin{equation}{section}
\newtheorem*{acknow}{Acknowledgements}

 \newcommand{\R}{\mathbb{R}}
 \newcommand{\Pol}{\mathrm{P}}
 \newcommand{\F}{\mathcal{F}}
 \newcommand{\To}{\longrightarrow}

 \newcommand{\Z}{\mathbb{Z}}

\begin{document}

\title[Cohomological rigidity and the number of homeomorphism types]
 {Cohomological rigidity and the number of homeomorphism types for small covers over prisms}

\author{Xiangyu Cao and  Zhi L\"u }

\address{School of Mathematical Sciences, Fudan University, Shanghai, P. R. China}
\email{xiangyu.cao08@gmail.com}
\address{School of Mathematical Sciences and The Key Laboratory of Mathematics for Nonlinear Sciences of Ministry of Education,  Fudan
 University, Shanghai, 200433, P.R. China.}
\email{zlu@fudan.edu.cn}
\thanks{This work is supported by grants from FDUROP (No. 080705),  NSFC (No. 10671034, No. J0730103 and  No. 10931005) and Shanghai NSF (No. 10ZR1403600)}

\begin{abstract}
In this paper, based upon the basic theory for glued manifolds in M.W. Hirsch (1976) \cite[Chapter 8, \S 2 Gluing Manifolds Together]{h}, we give a method of constructing homeomorphisms
between two small covers over simple convex polytopes. As a result
we classify, up to homeomorphism, all small covers over a
3-dimensional prism $\Pol^3(m)$ with $m\geq 3$. We introduce two
invariants from colored prisms and other two invariants from
ordinary cohomology rings with $\Z_2$-coefficients  of small covers. These invariants can form
a complete invariant system of homeomorphism types of all small
covers over a prism in most cases. Then we  show that the
cohomological rigidity holds for all small covers over a prism
$\Pol^3(m)$ (i.e., cohomology rings with $\Z_2$-coefficients of all small covers over a
$\Pol^3(m)$ determine their homeomorphism types). In addition, we
also calculate the number of homeomorphism types of all small covers
over  $\Pol^3(m)$.
\end{abstract}

\maketitle

\section{Introduction}
In 1991,  Davis and  Januszkiewicz~\cite{dj} introduced and studied
a class of $({\Bbb Z}_2)^n$-manifolds (called small covers), which
belong to the topological version of real toric varieties.
 An $n$-dimensional small cover $M^n$ is a  closed $n$-manifold  with a locally standard
$({\Bbb Z}_2)^n$-action such that its orbit space is a simple convex
$n$-polytope $P^n$. As shown in \cite{dj}, $P^n$ naturally admits a
characteristic function $\lambda$ defined on the facets of $P^n$
(here we also call $\lambda$ a $({\Bbb Z}_2)^n$-coloring on $P^n$), so
that the geometrical topology and the algebraic topology of the
small cover $M^n$ can be completely determined by the pair $(P^n,
\lambda)$. In other words, the Davis--Januszkiewicz theory for small
covers indicates the following two key points:
\begin{enumerate}
\item[$\bullet$]  Each small cover $\pi: M^n\longrightarrow P^n$  can be reconstructed from $(P^n,
\lambda)$, and this reconstruction is denoted by $M(\lambda)$. Thus,
all small covers over a simple convex polytope $P^n$ correspond to
all $({\Bbb Z}_2)^n$-colorings on $P^n$.
\item[$\bullet$]  The  algebraic topology of a small cover $\pi:
M^n\longrightarrow P^n$, such as equivariant cohomology, mod 2 Betti
numbers and ordinary cohomology,  etc., can be explicitly expressed
in terms of $(P^n, \lambda)$.
\end{enumerate}
In the recent years, much further research on small covers has been
carried on (see, e.g. \cite{i}, \cite{gs}, \cite{nn},  \cite{ccl},
\cite{c}, \cite{cmo}, \cite{lm}, \cite{ly}, \cite{km}, \cite{m1}--\cite{m2}). In
some sense, the classification up to equivariant homeomorphism of
small covers over a simple convex polytope has been understood very
well. Actually, this can be seen from the following two kinds of
viewpoints: One is that two small covers $M(\lambda_1)$ and
$M(\lambda_2)$ over a simple convex polytope $P^n$ are equivariantly
homeomorphic if and only if there is an automorphism $h\in \text{Aut}(P^n)$
such that $\lambda_1=\lambda_2\circ \bar{h}$ where $\bar{h}$ induced
by $h$ is an automorphism on all facets of $P^n$ (see \cite{lm});
the other one is that two small covers $M(\lambda_1)$ and $M(\lambda_2)$
over a simple convex polytope $P^n$ are equivariantly homeomorphic
if and only if their equivariant cohomologies are isomorphic as $H^*(B({\Bbb
Z}_2)^n; {\Bbb Z}_2)$-algebras (see \cite{m1}). However, in
non-equivariant case, the classification up to homeomorphism of
small covers over a simple convex polytope is far from understood
very well except for few special polytopes (see, e.g. \cite{gs},
\cite{km}, \cite{m2}, \cite{cmo}, \cite{ly}).

In this paper we shall introduce an approach (called the {\em sector
method}) of constructing homeomorphisms between two small covers.
The basic idea of sector method is simply stated as follows: we
first cut a $(\Z_2)^n$-colored simple convex $n$-polytope $(P^n,
\lambda)$ in ${\Bbb R}^n$ into two parts $P_1$ and $P_2$ by using an
$(n-1)$-dimensional hyperplane $H$ such that the section $S$ cut out
by $H$ is an $(n-1)$-dimensional simple convex polytope with certain
property. Of course, $S$ naturally inherits a coloring from $(P^n,
\lambda)$. Then by using automorphisms of $S$ and automorphisms of
$({\Bbb Z}_2)^n$ we can construct  new colored polytopes $(P'^n,
\lambda')$ from $(P^n, \lambda)$ (note that generally $P'^n$ may not
be combinatorially equivalent to $P^n$), and further obtain new
small covers $M(\lambda')$ from those new colored polytopes $(P'^n,
\lambda')$ by the reconstruction of small covers. Moreover, based upon the basic theory for glued manifolds in \cite[Chapter 8, \S 2 Gluing Manifolds Together]{h}, we
study how to construct the homeomorphisms between $M(\lambda)$ and
$M(\lambda')$. In particular, we give the method of constructing the
homeomorphisms between $M(\lambda)$ and $M(\lambda')$ (see
Theorems~\ref{extend}--\ref{t:GENERAL}).

As an application, up to homeomorphism we shall classify all small
covers over prisms. Let $\Pol^3(m)$ denote a 3-dimensional prism
that is the product of $[0,1]$ and an $m$-gon  where $m\geq 3$. Let
$\Lambda(\Pol^3(m))=\{\lambda\big| \lambda \text{ is a $(\Z_2)^3$-coloring on $\Pol^3(m)$}\}$, and let
$\Gamma(\Pol^3(m))=\{M(\lambda)\big|\lambda\in\Lambda(\Pol^3(m))\}$.
Using the sector method, we first study rectangular sectors  and
find seven rectangular sectors with a good twist
$\Psi(\psi,\rho,v_0)$ in Section~\ref{Rectangular} (see
Subsection~\ref{hom} for $\Psi(\psi,\rho,v_0)$). We then use
these seven rectangular sectors to define some operations on the
 coloring sequences of side-faces of colored polytopes $(\Pol^3(m),
\lambda)$ in Section~\ref{prism}. Furthermore, we show that all
$(\Pol^3(m), \lambda)$ can be reduced to some canonical forms
without changing the homeomorphism type of the small covers
$M(\lambda)$ (see Propositions~\ref{c1}--\ref{c3} and \ref{pNONTRI}).
In particular, we introduce two combinatorial invariants $m_\lambda$
and $n_\lambda$ in Section~\ref{prism}, and then show in
Section~\ref{cohomology} that $(m_\lambda, n_\lambda)$ is actually a
complete homeomorphism invariant of a class of small covers
$M(\lambda)$ (see Corollary~\ref{complete invariant}). In addition,
we also introduce two algebraic invariants $\Delta(\lambda)$ and
${\mathcal{B}}(\lambda)$ in $H^*(M(\lambda); \Z_2)$, which can
become a complete homeomorphism invariant in most cases.
$(m_\lambda, n_\lambda)$ and $(\Delta(\lambda),
{\mathcal{B}}(\lambda))$ are of interest because of the nature of
colored polytopes. With the help of these invariants, we can obtain
the following cohomological rigidity result.

\begin{thm}[Cohomological  rigidity]\label{RIGIDITY}
Two small covers $M(\lambda_1)$ and $M(\lambda_2)$ in
$\Gamma(\Pol^3(m))$ are homeomorphic if and only if  $H^*(M(\lambda_1);{\Bbb Z}_2)$ and
$H^*(M(\lambda_2);{\Bbb Z}_2)$ are isomorphic as graded rings.
\end{thm}
\begin{rem}
  Kamishima and  Masuda  showed in \cite{km} that the
cohomological rigidity holds for small covers over a cube. In
addition, Masuda in \cite{m2} also gave a cohomological non-rigidity
example of dimension 25. A further question is what extent the
cohomological rigidity can extend to small covers.
\end{rem}

In addition, we  also determine the number of homeomorphism classes
of all small covers in $\Gamma(\Pol^3(m))$.

\begin{thm} [Number of homeomorphism classes]\label{NUMBER}
Let $N(m)$ denote the number of homeomorphism classes of all small
covers in $\Gamma(\Pol^3(m))$. Then
$$N(m)=\begin{cases}
2 & \text{ if $m=3$}\\
4 & \text{ if $m=4$}\\
\sum_{0\leq k\leq
{m\over 2}}([{k\over 2}]+1)+6 & \text{ if $m>4$ is even}\\
\sum_{1\leq k\leq
{m\over 2}}([{k\over 2}]+1)+4 & \text{ if $m>4$ is odd.}\\
\end{cases}$$
\end{thm}

The arrangement of this paper is as follows. In Section~\ref{theory}
 we review the basic theory about small covers. In Section~\ref{sector} we
introduce the sector method. In Section~\ref{Rectangular} we apply
the sector method to colored prisms and discuss  the rectangular
sectors.  In Section~\ref{prism} we use the rectangular sectors to
define some operations on the coloring sequences of  colored prisms
and then determine the canonical forms of all colored prisms. In particular, we also introduce two combinatorial invariants $m_\lambda$
and $n_\lambda$ therein. In
Section~\ref{cohomology} we introduce two algebraic invariants of
cohomology and then calculate such two invariants of all small
covers. Finally, we complete the proofs of Theorems~\ref{RIGIDITY}
and~\ref{NUMBER} in Section~\ref{proof}.

\section{Theory of small covers}\label{theory}

The purpose of this section is to briefly review the theory of small
covers. Throughout the following assume that $\pi:
M^n\longrightarrow P^n$ is a small cover over a simple convex
$n$-polytope $P^n$. Note  that a simple convex $n$-polytope $P^n$
means that exactly $n$ faces of codimension-one (i.e., facets) meet
at each of its vertices. Let $\mathcal{F}(P^n)=\{F_1,...,F_\ell\}$
denote the set of all facets of $P^n$.

\subsection{Coloring and  reconstruction}
Take a $k$-face $F^k$ of $P^n$, an easy observation shows (see also
\cite[Lemma 1.3]{dj}) that $\pi^{-1}(F^k)\longrightarrow F^k$ is
still a $k$-dimensional small cover. In particular, for any $x\in
\pi^{-1}(\text{int}F^k)$, its isotropy subgroup $G_x$ is independent
of the choice of $x$, denoted by $G_F$.  $G_F$ is isomorphic to
$(\Z_2)^{n-k}$, and $G_F$ fixes $\pi^{-1}(F^k)$ in $M^{n}$. In the
case $k=n-1$, $F^{n-1}$ is a facet and $G_F$ has rank $1$, so that
 $G_F$ uniquely corresponds to a nonzero
vector $v_F$ in $({\Bbb Z}_2)^n$. Then there is a natural map (called
{\em characteristic function})
$$\lambda:\mathcal{F}(P)\longrightarrow (\Z_2)^n$$ by mapping each facet $F$ to its corresponding nonzero vector $v_F$ in $({\Bbb Z}_2)^n$  with the property $(\star)$: whenever the intersection of
some facets $F_{i_1}, ..., F_{i_r}$ in $\mathcal{F}(P^n)$ is
nonempty, $\lambda(F_{i_1}),..., \lambda(F_{i_r})$ are linearly
independent in $(\Z_2)^n$. Note that if  each nonzero vector of
$(\Z_2)^n$ is regarded as being a color, then the characteristic
function $\lambda$ means that each facet is colored by a color.
Thus, we also call $\lambda$ a {\em $(\Z_2)^n$-coloring} on $P^n$
here. By $\Lambda(P^n)$ we denote the set of all $(\Z_2)^n$-colorings
on $P^n$.
\begin{rem}
Since $P^n$ is simple, for each $k$-face $F^k$, there are $n-k$
facets $F_{i_1},...,F_{i_{n-k}}$ such that
$F^k=F_{i_1}\cap\cdots\cap F_{i_{n-k}}$ and $\pi^{-1}(F^k)$ is a
transversal intersection of $\pi^{-1}(F_{i_1})$, $...,$
$\pi^{-1}(F_{i_{n-k}})$. Then the group $G_F$ determined by $F^k$ is
actually generated by $\lambda(F_{i_1})$, $...,$
$\lambda(F_{i_{n-k}})$.
\end{rem}
Davis and Januszkiewicz \cite{dj} gave a reconstruction  of $M^n$ by
using the $(\Z_2)^n$-coloring $\lambda$ and the product bundle
$P^n\times (\Z_2)^n$ over $P^n$. Geometrically this reconstruction is
exactly formed by gluing $2^n$ copies of $P^n$ along their boundaries
via $\lambda$. This reconstruction can be written as follows:
$$M(\lambda):=P^n\times (\Z_2)^n/(p, v)\sim (p, v+\lambda(F)) \text{ for } p\in F\in \mathcal{F}(P^n).$$
Then we have
\begin{thm}[Davis--Januszkiewicz] All small covers over
$P^n$ are given by $\Gamma(P^n)$$ =\{M(\lambda)|\lambda\in
\Lambda(P^n)\}$.
\end{thm}
There is a natural action of $\text{GL}(n,\Z_2)$ on $\Lambda(P^n)$
defined by $\lambda\longmapsto \alpha\circ \lambda$, and it is easy
to see that such an action is free, and it also induces an action of
$\text{GL}(n,\Z_2)$ on $\Gamma(P^n)$ by $M(\lambda)\longmapsto
M(\alpha\circ \lambda)$. Following \cite{dj}, two small covers
$M(\lambda_1)$ and $M(\lambda_2)$ in $\Gamma(P^n)$ are said to be
{\em Davis--Januszkiewicz equivalent} if there is an $\alpha\in
\text{GL}(n,\Z_2)$ such that $\lambda_1=\alpha\circ \lambda_2$.
Thus, each Davis--Januszkiewicz equivalence class in $\Gamma(P^n)$ is
actually an orbit of the action of $\text{GL}(n,\Z_2)$ on
$\Gamma(P^n)$.

\subsection{Betti numbers and $h$-vector}
The notion of the  $h$-vector plays an  essential important role in
the  theory of polytopes, while the notion of Betti numbers is also
so
 important
in the topology of manifolds. Davis--Januszkiewicz theory indicates
that the Dehn--Sommerville relations for the $h$-vectors and the
Poincar\'e duality for the mod 2 Betti numbers are essentially consistent
in the setting of small covers.

Let  $P^*$ be the dual of $P^n$ that is a simplicial polytope. Then
the boundary $\partial P^*$ is  a finite simplicial
complex of dimension $n-1$. For $0\leq i\leq n-1$, by $f_i$ we
denote the number of all $i$-faces in $\partial P^*$. Then the vector $(f_0,
f_1,...,f_{n-1})$ is called the {\em $f$-vector} of $P^n$, denoted
by ${\bf f}(P^n)$. Then the {\em $h$-vector} denoted by ${\bf
h}(P^n)$ of $P^n$ is an integer vector $(h_0, h_1,..., h_n)$ defined
from the following equation $$ h_0t^n + \cdots + h_{n-1}t + h_n =
(t-1)^n + f_0(t-1)^{n-1}+\cdots + f_{n-1}.$$
\begin{thm}[Davis--Januszkiewicz]\label{betti}
Let $\pi: M^{n}\longrightarrow P^n$ be a small cover over a simple
convex polytope $P^n$. Then
 $${\bf h}(P^n)=(h_0,...,h_n)=(b_0,...,b_n)$$
where $b_i=\dim H^i(M^n;{\Bbb Z}_2)$.
\end{thm}
\begin{rem}
We see from Theorem~\ref{betti} that the Poincar\'e duality
$b_i=b_{n-i}$ agrees with the Dehn-Sommerville relation
$h_i=h_{n-i}$.
\end{rem}
\begin{example}\label{betti}
For the prism $\Pol^3(m)$, ${\bf h}(\Pol^3(m))=(1, m-1,  m-1, 1)$,
so any small cover over $\Pol^3(m)$ has the mod 2 Betti numbers
$(b_0, b_1, b_2, b_3)=(1, m-1,  m-1, 1)$.
\end{example}

\subsection{Stanley--Reisner face ring and equivariant cohomology}
Stanley--Reisner face ring is a basic combinatorial invariant, and
equivariant cohomology is  an essential invariant in the theory of
transformation groups. Davis--Januszkiewicz theory  indicates that
these two kinds of invariants are also essentially consistent in the
setting of small covers.

Let $P^n$ be a simple convex polytope with facet set
$\mathcal{F}(P^n)=\{F_1,...,F_\ell\}$. Following \cite{dj},  the
Stanley--Reisner face ring of $P^n$ over $\Z_2$, denoted by
$\Z_2(P^n)$, is defined as follows:
$$\Z_2(P^n)=\Z_2[F_1,...,F_\ell]/I$$
where the $F_i$'s are regarded as indeterminates of degree one, and
$I$ is a homogeneous ideal generated by all sequence free monomials
of the form $F_{i_1}\cdots F_{i_s}$ with $F_{i_1}\cap\cdots\cap
F_{i_s}=\emptyset$.

 \begin{example}
Let $P^n$ be an $n$-simplex $\Delta^n$ with $n+1$ facets
$F_1,...,F_{n+1}$. Then
$$\Z_2(\Delta^n)=\Z_2[F_1,...,F_{n+1}]/(F_1\cdots F_{n+1}).$$
 \end{example}

 \begin{example}
Let $F_1, ...,  F_{2n}$ be $2n$ facets of an $n$-cube $I^n$ with
$F_i\cap F_{i+n}=\emptyset, i=1,..., n$. Then
$\Z_2(I^n)=\Z_2[F_1, ..., F_{2n}]/(F_iF_{i+n}|i=1,...,n).$
 \end{example}

\begin{thm}[Davis--Januszkiewicz]\label{d-j}
Let $\pi: M^{n}\longrightarrow P^n$ be a small cover over a simple
convex polytope $P^n$. Then its equivariant cohomology is
$$H^*_{(\Z_2)^n}(M^{n};\Z_2)\cong \Z_2(P^n).$$
\end{thm}

\subsection{Ordinary cohomology} Let $\pi: M^{n}\longrightarrow P^n$ be a small cover over a simple
convex polytope $P^n$ with $\mathcal{F}(P^n)=\{F_1,...,F_\ell\}$,
and $\lambda: \mathcal{F}(P^n)\longrightarrow (\Z_2)^n$ its
$(\Z_2)^n$-coloring.  Now let us extend $\lambda:
\mathcal{F}(P^n)\longrightarrow (\Z_2)^n$ to a linear map
$\widetilde{\lambda}: (\Z_2)^\ell\longrightarrow (\Z_2)^n$ by replacing
$\{F_1,...,F_\ell\}$ by the basis $\{e_1,...,e_\ell\}$ of
$(\Z_2)^\ell$. Then $\widetilde{\lambda}: (\Z_2)^\ell\longrightarrow
(\Z_2)^n$ is surjective, and $\widetilde{\lambda}$ can be regarded as
an $n\times \ell$-matrix $(\lambda_{ij})$, which is written as
$(\lambda(F_1), ...,\lambda(F_\ell)).$ It is well-known that
$H_1(B(\Z_2)^\ell;\Z_2)=H_1(E(\Z_2)^n\times_{(\Z_2)^n}
M^n;\Z_2)=(\Z_2)^\ell$ and $H_1(B(\Z_2)^n;\Z_2)=(\Z_2)^n$. So one has that
$p_*: H_1(E(\Z_2)^n\times_{(\Z_2)^n} M^n;\Z_2)\longrightarrow
H_1(B(\Z_2)^n;\Z_2)$ can be identified with $\widetilde{\lambda}:
\Z_2^\ell\longrightarrow (\Z_2)^n$, where $p:E(\Z_2)^n\times_{(\Z_2)^n}
M^n\longrightarrow B(\Z_2)^n$ is the fibration of the Borel
construction associating to the universal principal $(\Z_2)^n$-bundle
$E(\Z_2)^n\longrightarrow B(\Z_2)^n$. Furthermore,  $p^*:
H^1(B(\Z_2)^n;\Z_2)\longrightarrow H^1_{(\Z_2)^n}(M^n;\Z_2)$ is
identified with the dual map $\widetilde{\lambda}^*:
(\Z_2)^{n*}\longrightarrow (\Z_2)^{\ell*}$, where
$\widetilde{\lambda}^*=\widetilde{\lambda}^\top$ as  matrices.
Therefore, column vectors of $\widetilde{\lambda}^*$ can be
understood as linear combinations of $F_1,...,F_\ell$ in the face
ring $\Z_2(P^n)=\Z_2[F_1,...,F_\ell]/I$. Write
$\lambda_i=\lambda_{i1}F_1+\cdots+\lambda_{i\ell}F_\ell.$
Let $J_\lambda$ be the homogeneous ideal $(\lambda_1,...,\lambda_n)$
in $\Z_2[F_1,...,F_\ell]$. Davis and Januszkiewicz calculated the
ordinary cohomology of $M^n$, which is stated as follows.

\begin{thm}[Davis--Januszkiewicz]
Let $\pi: M^{n}\longrightarrow P^n$ be a small cover over a simple
convex polytope $P^n$. Then its ordinary cohomology
$$H^*(M^{n};\Z_2)\cong \Z_2[F_1,...,F_\ell]/I+J_\lambda.$$
\end{thm}

The following result which will be used later is due to Nakayama and
Nishimura \cite{nn}.

\begin{prop}\label{orient}
Let $\pi: M^{n}\longrightarrow P^n$ be a small cover over a simple
convex polytope $P^n$, and $\lambda:
\mathcal{F}(P^n)=\{F_1,...,F_\ell\}\longrightarrow(\Z_2)^n$ its
$(\Z_2)^n$-coloring. Then $M^n$ is orientable if and only if there
exists an automorphism $\sigma\in \text{\rm GL}(n,\Z_2)$ such that
$\lambda'=\sigma\circ \lambda$ satisfies $\sum_{j=1}^n
\lambda'_{jl}\equiv 1\mod 2$ for all $1\leq l\leq \ell$, where
$\widetilde{\lambda'}=(\lambda'_{ij}):(\Z_2)^\ell\longrightarrow
(\Z_2)^n$ is the linear extension of $\lambda'$, as  above.
\end{prop}

\section{Sector Method}\label{sector}

Each point of a simple convex polytope $P^n$ has a neighborhood
which is affinely isomorphic to an open subset of ${\Bbb R}^n_{\geq
0}$, so $P^n$ is an $n$-dimensional nice manifold with corners (see
\cite{d}). An automorphism of $P^n$ is a self-homeomorphism of $P^n$
as a manifold with corners, and by $\text{Aut}(P^n)$ we denote the
group of automorphisms of $P^n$. All faces of $P^n$ form a poset by
inclusion. An automorphism of $\mathcal{F}(P^n)$ is a bijection from
$\mathcal{F}(P^n)$ to itself which preserves the poset structure of
all faces of $P^n$, and by $\text{Aut}(\mathcal{F}(P^n))$ we denote
the group of automorphisms of $\mathcal{F}(P^n)$. Each automorphism
of $\text{Aut}(P^n)$ naturally induces an automorphism of
$\mathcal{F}(P^n)$. It is well-known (see \cite{bp} or \cite{z})
that two simple convex polytopes are combinatorially equivalent  if
and only if they are homeomorphic as manifolds with corners. Thus,
the natural homomorphism $\Phi: \text{Aut}(P^n)\longrightarrow
\text{Aut}(\mathcal{F}(P^n))$ is surjective. Throughout the following, we use the convention that all simple $n$-polytopes are naturally embedded in ${\Bbb R}^n$.

\begin{defn} Let $P^n$ be a simple $n$-polytope and $S$ be a simple
$(n-1)$-polytope. An embedding $i:S\To P^n$ is called a \emph{sector}
if the following conditions are satisfied:

\begin{enumerate}
\item[(a)] $i(S)$ divides $P^n$ into two simple $n$-polytopes $P_1$ and $P_2$
 such that $i(S)$ is the common facet of $P_1$ and $P_2$.
\item[(b)]For every $k$-face $F\subset S$, there is a unique $(k+1)$-dimensional face $H$ of $P^n$ such that
$i(F)\subset H$, where the uniqueness means that $S$ contains no any vertex of $P$.
\end{enumerate}
\end{defn}

A sector $i:S\To P^n$ naturally induces a poset structure-preserving map $i_*:\F(S)\rightarrow\F(P)$. If
$P^n$ admits a $(\Z_2)^n$-coloring $\lambda$, then $\bar{\lambda}=\lambda\circ i_*$ gives a $(\Z_2)^n$-coloring on $S$, called the {\em derived
coloring}. Obviously,  $\bar{\lambda}$  satisfies the independence condition:
whenever the intersection of some facets $F_{l_1}, ..., F_{l_r}$ in
$\mathcal{F}(S)$ is nonempty, $\bar{\lambda}(F_{l_1}), ...,
\bar{\lambda}(F_{l_r})$ are linearly independent in $(\Z_2)^n$. Let $(\psi, \rho)$ be a pair of  $\psi \in \text{Aut}(S)$ and $\rho\in \text{GL}(n,
\Z_2)$. If
$\bar{\lambda}\circ\psi=\rho\circ\bar{\lambda}$, then $(\psi,\rho)$ is called an \emph{isomorphism pair} of preserving $(S, \bar{\lambda})$, where we just abuse
$\psi$ and the  automorphism of $\mathcal{F}(S)$ induced by
$\psi$.

\subsection{Construction of new colored simple polytopes}
 Given a $(\Z_2)^n$-colored simple $n$-polytope $(P,\lambda)$, let $i:S\To P$ be a sector with $P_1$ and $P_2$  two
polytopes cut out by $i(S)$, and let $i_r:S\rightarrow P_r, r=1,2$, be the embeddings
induced by $p\mapsto i(p)$. Then $P=P_1\coprod P_2/i_1(s)\sim
i_2(s).$

Now choose  an isomorphism pair $(\psi,\rho)$ of preserving $(S, \bar{\lambda})$.
 Suppose that
$P'$ is another simple $n$-polytope with $j: S\rightarrow P'$
another sector, cutting $P'$ into ${P'}_1$ and ${P'}_2$, such that
there are $g_r: P_r\rightarrow{P'}_r, r=1,2$,  which are
combinatorially equivalent with $g_1\circ i_1=j_1$ and $g_2\circ
i_2\circ\psi=j_2$. Thus
\begin{equation}\label{constr}
P'=P_1'\coprod P_2'/j_1(s)\sim
j_2(s)=g_1(P_1)\coprod g_2(P_2)/g_1\circ i_1(s)\sim g_2\circ
i_2\circ \psi(s).
\end{equation}
It is not difficult to see that generally $P$ is not combinatorially equivalent to $P'$ although
$P_r$ is combinatorially equivalent to $P'_r$ $(r=1,2)$. In particular, if $g_r$ is chosen as the identity or a mirror reflection along an $(n-1)$-dimensional hyperplane in ${\Bbb R}^n$, then
$P'$ is exactly constructed by gluing $P_1$ and $P_2$ along their facets $i_1(S)$ and $i_2(\psi(S))$ via $\psi$.

Next, we shall give a $(\Z_2)^n$-coloring on $P'$ via $\lambda$ and $\rho$.
For each facet $F\in\F(P')$, if $F\cap{P_1}'\neq\emptyset$,
then there is a unique facet $F_1\in \mathcal{F}(P)$ such that
$g_1(F_1\cap P_1)\subset F$. Similarly, if
$F\cap{P'}_2\neq\emptyset$, then there is also a unique facet
$F_2\in \mathcal{F}(P)$ such that $g_2(F_2\cap P_2)\subset F$.
 Furthermore,   define
 $\lambda': \F(P')\rightarrow(\Z_2)^n$ in the following way:
\begin{equation}\label{coloring}
\lambda'(F)=\begin{cases}
\lambda(F_1) & \text{ if $F\cap{P_1}'\neq\emptyset$}\\
\rho^{-1}\circ\lambda(F_2) & \text{ if $F\cap{P_2}'\neq\emptyset$.}
\end{cases}
\end{equation}
 Such $\lambda'$ is well-defined. In fact,
if $F$ has nonempty intersection with both ${P'}_1$ and $P'_2$, then
$F$ must lie in the image of $j_*$, say $F=j_*(F_s)$ where $F_s\in
\mathcal{F}(S)$.  Then from $g_1\circ i_1=j_1$ we see that
$F_1=i_*(F_s)$, and from $g_2\circ i_2\circ\psi=j_2$ we see that
$F_2=i_*(\psi(F_s))$. So
$\lambda(F_1)=\bar{\lambda}(F_s)=\rho^{-1}\circ\bar{\lambda}\circ\psi(F_s)=\rho^{-1}\circ\lambda\circ
i_*\circ\psi(F_s)=\rho^{-1}\circ\lambda(F_2)$ as desired.

\vskip .2cm In summary,
we now have two colored polytopes $(P,\lambda)$ and
$({P}',\lambda')$. Then using the reconstruction method of small covers, we
obtain two small covers $M(\lambda)=P\times (\Z_2)^n/(p, v)\sim (p,
v+\lambda(F)) \text{ for } p\in F\in \mathcal{F}(P)$ and
$M(\lambda')=P'\times (\Z_2)^n/(p, v)\sim (p, v+\lambda'(F)) \text{
for } p\in F\in \mathcal{F}(P')$. In addition, we have also a colored polytope $(S, \bar{\lambda})$. We can still use  the reconstruction technique  of small covers to $(S, \bar{\lambda})$ to get $$\mathcal{S}=S\times(\Z_2)^n/\{(s,v)\sim(s,v+\bar{\lambda}(f))| s\in f\in\F(S)\}.$$
Then it is easy to see that $\mathcal{S}$ is an $(n-1)$-dimensional
closed manifold (but possibly disconnected), called the {\em sector
manifold} here. As noted in \cite[Remark 5.1]{ly} for the 2-dimensional case,  $\mathcal{S}$ is not a small cover in the sense of Davis-Januszkiewicz.

\subsection{Relation between homeomorphism types of $M(\lambda)$ and $M(\lambda')$}
 \label{hom} Following the above notions, we shall study when $M(\lambda)$ and $M(\lambda')$ are
homeomorphic by using the basic theory for glued manifolds of \cite[Chapter 8, \S 2 Gluing Manifolds Together]{h}.

 Let $\pi: M(\lambda)\longrightarrow P$ and $\pi' : M(\lambda') \longrightarrow P'$ be the projections.  Set $M_r={\pi}^{-1}(P_r), r=1,2$. Then each $M_r$ is a manifold with boundary $\pi^{-1}(i(S))$.
  The sector $i: S\longrightarrow P$ naturally induces
an embedding $\iota: \mathcal{S}\hookrightarrow M(\lambda)$ defined
by $x=\{(s,v)\}\longmapsto \{(i(s), v)\}$, and $i_r: S\longrightarrow
P_r$ also induces the embedding $\iota_r: \mathcal{S}\hookrightarrow
M_r$ $(r=1,2)$. Obviously, $\partial
M_r=\iota_r(\mathcal{S})={\pi}^{-1}(i(S))=\iota(\mathcal{S})$. Using
this terminology one can write
$$M(\lambda)=M_1\coprod
M_2/\iota_1(x)\sim\iota_2(x)=M_1\coprod M_2/\{(i_1(s), v)\}\sim
\{(i_2(s), v)\},$$ i.e., $M(\lambda)$ is obtained by gluing  $M_1$
and $M_2$ together along their common boundary $\iota(\mathcal{S})$
via $\iota_1$ and $\iota_2$.
 Similarly, set $M'_r={\pi'}^{-1}(P'_r), r=1,2$. We have  the embeddings
$\iota': \mathcal{S}\hookrightarrow M(\lambda')$ induced by
$j:S\longrightarrow P'$ and  $\iota'_r: \mathcal{S}\hookrightarrow M'_r$ induced by
$j_r:S\longrightarrow P'_r$ $(r=1,2)$ such that $\partial
M'_r=\iota'_r(\mathcal{S})={\pi'}^{-1}(j(S))=\iota'(\mathcal{S})$.  Then one has
that $$ M(\lambda') = M'_1\coprod M'_2/\iota'_1(x)\sim \iota'_2(x) =
M'_1\coprod M'_2/\{(j_1(s), v)\}\sim \{(j_2(s), v)\}.
$$
By (\ref{constr}) and (\ref{coloring}), we see that
 $g_1:P_1\longrightarrow
P'_1$ may induce an equivariant homeomorphism $\widetilde{g}_1:
M_1\longrightarrow M'_1$ by mapping $\{(p,v)\}\longmapsto
\{(g_1(p), v)\}$, while  $g_2:P_2\longrightarrow P'_2$ may induce a
weakly equivariant homeomorphism $\widetilde{g}_2:
M_2\longrightarrow M'_2$ by mapping $\{(p,v)\}\longmapsto
\{(g_2(p), \rho^{-1}(v))\}$.

On the other hand, take any element  $v_0\in (\Z_2)^n$, using the
relation $\bar{\lambda}\circ\psi=\rho\circ\bar{\lambda}$ we see
easily that  $(\psi, \rho)$ also induces
a weakly equivariant homeomorphism $\Psi:
\mathcal{S}\rightarrow\mathcal{S}$ defined by
$\{(s,v)\}\mapsto\{(\psi(s),\rho(v)+v_0)\}$, so that  $\iota_2\circ
\Psi$ gives a new embedding of $\mathcal{S}$ in $M_2$ and
$\widetilde{g}_2\circ \iota_2\circ \Psi$ also gives a new embedding
of $\mathcal{S}$ in $M'_2$.  Let $\widetilde{M}(\lambda)=M_1\coprod
M_2/\iota_1(x)\sim\iota_2(\Psi(x))$. Then one has that

\begin{lem} \label{induced}
 $\widetilde{M}(\lambda)$
is homeomorphic to $M(\lambda')$.
\end{lem}

\begin{proof} Let $z=\{(p, v)\}\in \widetilde{M}(\lambda)$.
Define $\Pi: \widetilde{M}(\lambda)\longrightarrow M(\lambda')$ by
$$\Pi(z)=\begin{cases}
\widetilde{g}_1(z) & \text{if $z\in M_1$}\\
\widetilde{g'}_2(z) & \text{if $z\in M_2$}
\end{cases}$$
where $\widetilde{g'}_2(z)=\{(g_2(p),
\rho^{-1}(v)+\rho^{-1}(v_0))\}$. To show that $\Pi$ is a
homeomorphism,   it suffices to prove that for $x\in \mathcal{S}$, if
$\iota_1(x)\sim\iota_2(\Psi(x))$ in $\widetilde{M}(\lambda)$, then
$\Pi\circ \iota_1(x)\sim\Pi\circ\iota_2(\Psi(x))$ in $M(\lambda')$.
Since $\Pi\circ \iota_1=\widetilde{g}_1\circ \iota_1=\iota'_1$, it
needs only to check that $\Pi\circ\iota_2\circ\Psi=\iota'_2$. Let
$x=\{(s,v)\}$. Then
$\Pi\circ\iota_2\circ\Psi(x)=\widetilde{g'}_2\circ \iota_2\circ
\Psi(x)=\widetilde{g'}_2\circ \iota_2(\{(\psi(s),
\rho(v)+v_0)\})=\widetilde{g'}_2(\{(i_2\circ\psi(s),
\rho(v)+v_0)\}=\{(g_2\circ i_2\circ
\psi(s),\rho^{-1}(\rho(v)+v_0)+\rho^{-1}(v_0))\}=\{(g_2\circ
i_2\circ \psi(s),v)\}=\iota'_2(x)$, as desired.
\end{proof}

\begin{example}\label{special example}
Let $P$ be a $(\Z_2)^3$-colored 3-cube,   $P'$  a truncated 3-cube, and $S$ a square. Let $i: S\hookrightarrow P$ and
$j: S\hookrightarrow P'$ be two sectors, as shown in the following picture:
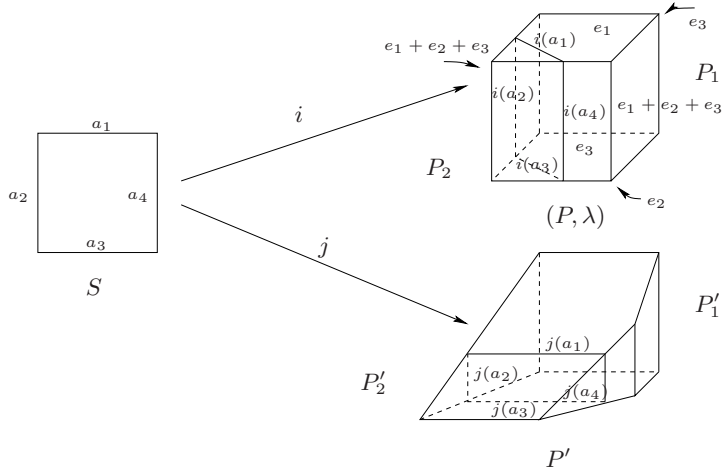
\[   \begin{picture}(0,0)%
\includegraphics{2.pstex}%
\end{picture}%
\setlength{\unitlength}{1973sp}%
\begingroup\makeatletter\ifx\SetFigFont\undefined%
\gdef\SetFigFont#1#2#3#4#5{%
  \reset@font\fontsize{#1}{#2pt}%
  \fontfamily{#3}\fontseries{#4}\fontshape{#5}%
  \selectfont}%
\fi\endgroup%
\begin{picture}(9394,5914)(1126,-5233)
\put(5551,-4186){\makebox(0,0)[lb]{\smash{{\SetFigFont{9}{10.8}{\rmdefault}{\mddefault}{\updefault}$P'_2$}}}}
\put(7876,-5161){\makebox(0,0)[lb]{\smash{{\SetFigFont{9}{10.8}{\rmdefault}{\mddefault}{\updefault}$P'$}}}}
\put(4726,-811){\makebox(0,0)[lb]{\smash{{\SetFigFont{9}{10.8}{\rmdefault}{\mddefault}{\updefault}$i$}}}}
\put(5026,-2461){\makebox(0,0)[lb]{\smash{{\SetFigFont{9}{10.8}{\rmdefault}{\mddefault}{\updefault}$j$}}}}
\put(7726,164){\makebox(0,0)[lb]{\smash{{\SetFigFont{6}{7.2}{\rmdefault}{\mddefault}{\updefault}$i(a_1)$}}}}
\put(7201,-511){\makebox(0,0)[lb]{\smash{{\SetFigFont{6}{7.2}{\rmdefault}{\mddefault}{\updefault}$i(a_2)$}}}}
\put(7876,-3661){\makebox(0,0)[lb]{\smash{{\SetFigFont{6}{7.2}{\rmdefault}{\mddefault}{\updefault}$j(a_1)$}}}}
\put(6976,-4036){\makebox(0,0)[lb]{\smash{{\SetFigFont{6}{7.2}{\rmdefault}{\mddefault}{\updefault}$j(a_2)$}}}}
\put(7201,-4486){\makebox(0,0)[lb]{\smash{{\SetFigFont{6}{7.2}{\rmdefault}{\mddefault}{\updefault}$j(a_3)$}}}}
\put(2101,-2386){\makebox(0,0)[lb]{\smash{{\SetFigFont{6}{7.2}{\rmdefault}{\mddefault}{\updefault}$a_3$}}}}
\put(2176,-886){\makebox(0,0)[lb]{\smash{{\SetFigFont{6}{7.2}{\rmdefault}{\mddefault}{\updefault}$a_1$}}}}
\put(1126,-1786){\makebox(0,0)[lb]{\smash{{\SetFigFont{6}{7.2}{\rmdefault}{\mddefault}{\updefault}$a_2$}}}}
\put(2626,-1786){\makebox(0,0)[lb]{\smash{{\SetFigFont{6}{7.2}{\rmdefault}{\mddefault}{\updefault}$a_4$}}}}
\put(2101,-2986){\makebox(0,0)[lb]{\smash{{\SetFigFont{9}{10.8}{\rmdefault}{\mddefault}{\updefault}$S$}}}}
\put(8476,314){\makebox(0,0)[lb]{\smash{{\SetFigFont{6}{7.2}{\rmdefault}{\mddefault}{\updefault}$e_1$}}}}
\put(8251,-1186){\makebox(0,0)[lb]{\smash{{\SetFigFont{6}{7.2}{\rmdefault}{\mddefault}{\updefault}$e_3$}}}}
\put(9676,389){\makebox(0,0)[lb]{\smash{{\SetFigFont{6}{7.2}{\rmdefault}{\mddefault}{\updefault}$e_3$}}}}
\put(9151,-1861){\makebox(0,0)[lb]{\smash{{\SetFigFont{6}{7.2}{\rmdefault}{\mddefault}{\updefault}$e_2$}}}}
\put(8776,-661){\makebox(0,0)[lb]{\smash{{\SetFigFont{6}{7.2}{\rmdefault}{\mddefault}{\updefault}$e_1+e_2+e_3$}}}}
\put(5851, 89){\makebox(0,0)[lb]{\smash{{\SetFigFont{6}{7.2}{\rmdefault}{\mddefault}{\updefault}$e_1+e_2+e_3$}}}}
\put(7876,-2086){\makebox(0,0)[lb]{\smash{{\SetFigFont{9}{10.8}{\rmdefault}{\mddefault}{\updefault}$(P, \lambda)$}}}}
\put(8101,-736){\makebox(0,0)[lb]{\smash{{\SetFigFont{6}{7.2}{\rmdefault}{\mddefault}{\updefault}$i(a_4)$}}}}
\put(8101,-4261){\makebox(0,0)[lb]{\smash{{\SetFigFont{6}{7.2}{\rmdefault}{\mddefault}{\updefault}$j(a_4)$}}}}
\put(7501,-1411){\makebox(0,0)[lb]{\smash{{\SetFigFont{6}{7.2}{\rmdefault}{\mddefault}{\updefault}$i(a_3)$}}}}
\put(6376,-1486){\makebox(0,0)[lb]{\smash{{\SetFigFont{9}{10.8}{\rmdefault}{\mddefault}{\updefault}$P_2$}}}}
\put(9751,-286){\makebox(0,0)[lb]{\smash{{\SetFigFont{9}{10.8}{\rmdefault}{\mddefault}{\updefault}$P_1$}}}}
\put(9751,-3211){\makebox(0,0)[lb]{\smash{{\SetFigFont{9}{10.8}{\rmdefault}{\mddefault}{\updefault}$P'_1$}}}}
\end{picture}%
\centering
   \]
where $\{e_1, e_2, e_3\}$ is the standard basis of $(\Z_2)^3$.
The $(\Z_2)^3$-coloring $\lambda$ on $P$ induces a $(\Z_2)^3$-coloring
$\bar{\lambda}$ on $S$ via $i$ as follows: $\bar{\lambda}:\mathcal{F}(S)=\{a_1,a_2,a_3,a_4\}\longrightarrow (\Z_2)^3$
by mapping $a_1\longmapsto e_1, a_2\longmapsto e_1+e_2+e_3, a_3\longmapsto e_2, a_4\longmapsto e_3$.  Choose an isomorphism pair $(\psi, \rho)$ of $(S, \bar{\lambda})$ as follows: $\psi:\mathcal{F}(S)\longrightarrow\mathcal{F}(S)$
is defined by $\psi(a_1)=a_2, \psi(a_2)=a_3, \psi(a_3)=a_4, \psi(a_4)=a_1$, and $\rho=\begin{bmatrix}
1 & 0 & 0\\
1 & 0 & 1\\
1 & 1 & 0
\end{bmatrix}$. Then we can have two combinatorial equivalences
$g_1:P_1\longrightarrow P'_1$ and $g_2: P_2\longrightarrow P'_2$ such that $g_1$ maps $i_1(a_1), i_1(a_2), i_1(a_3), i_1(a_4)$ into  $j_1(a_1), j_1(a_2), j_1(a_3), j_1(a_4)$ respectively, and $g_2$ maps $i_2(a_1),$ $i_2(a_2),$  $i_2(a_3),$  $i_2(a_4)$ into  $j_2(a_4), j_2(a_1), j_2(a_2), j_2(a_3)$ respectively. Furthermore, we obtain a $(\Z_2)^3$-coloring $\lambda'$ on $P'$,
as shown in the following picture:
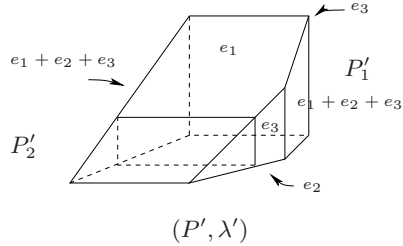
\[   \begin{picture}(0,0)%
\includegraphics{3.pstex}%
\end{picture}%
\setlength{\unitlength}{1973sp}%
\begingroup\makeatletter\ifx\SetFigFont\undefined%
\gdef\SetFigFont#1#2#3#4#5{%
  \reset@font\fontsize{#1}{#2pt}%
  \fontfamily{#3}\fontseries{#4}\fontshape{#5}%
  \selectfont}%
\fi\endgroup%
\begin{picture}(4969,3087)(5551,-5317)
\put(9751,-3211){\makebox(0,0)[lb]{\smash{{\SetFigFont{9}{10.8}{\rmdefault}{\mddefault}{\updefault}$P'_1$}}}}
\put(8176,-2911){\makebox(0,0)[lb]{\smash{{\SetFigFont{6}{7.2}{\rmdefault}{\mddefault}{\updefault}$e_1$}}}}
\put(5551,-3061){\makebox(0,0)[lb]{\smash{{\SetFigFont{6}{7.2}{\rmdefault}{\mddefault}{\updefault}$e_1+e_2+e_3$}}}}
\put(9151,-3586){\makebox(0,0)[lb]{\smash{{\SetFigFont{6}{7.2}{\rmdefault}{\mddefault}{\updefault}$e_1+e_2+e_3$}}}}
\put(9826,-2386){\makebox(0,0)[lb]{\smash{{\SetFigFont{6}{7.2}{\rmdefault}{\mddefault}{\updefault}$e_3$}}}}
\put(9226,-4636){\makebox(0,0)[lb]{\smash{{\SetFigFont{6}{7.2}{\rmdefault}{\mddefault}{\updefault}$e_2$}}}}
\put(7576,-5236){\makebox(0,0)[lb]{\smash{{\SetFigFont{9}{10.8}{\rmdefault}{\mddefault}{\updefault}$(P', \lambda')$}}}}
\put(8701,-3886){\makebox(0,0)[lb]{\smash{{\SetFigFont{6}{7.2}{\rmdefault}{\mddefault}{\updefault}$e_3$}}}}
\put(5551,-4186){\makebox(0,0)[lb]{\smash{{\SetFigFont{9}{10.8}{\rmdefault}{\mddefault}{\updefault}$P'_2$}}}}
\end{picture}%
\centering
   \]
   We know from \cite[Lemma 5.1]{ly} that the sector manifold $\mathcal{S}$ is a 2-dimensional torus $T^2$, so $M_1, M_2, M'_1, M'_2$ are all 2-dimensional solid  tori.
   In this case, we easily check that $\widetilde{M}(\lambda)$ and $M(\lambda')$ are not only homeomorphic, but also equivariantly homeomorphic.
\end{example}

It should be pointed out that Example~\ref{special example} only provides a very special case. Generally, $\widetilde{M}(\lambda)$ and $M(\lambda')$ may not be (weakly) equivariantly homeomorphic. Actually, for our purpose, we need to forget equivariant information on $M_1, M_2, M'_1, M'_2$ and $\mathcal{S}$. Let us look at the following example.

\begin{example}
Let $P$ and $P'$ be two prisms such that $P$ admits a $(\Z_2)^3$-coloring $\lambda$. Let $S$ be a triangle, which is embedded into $P$ and $P'$ respectively, as shown in the following picture:
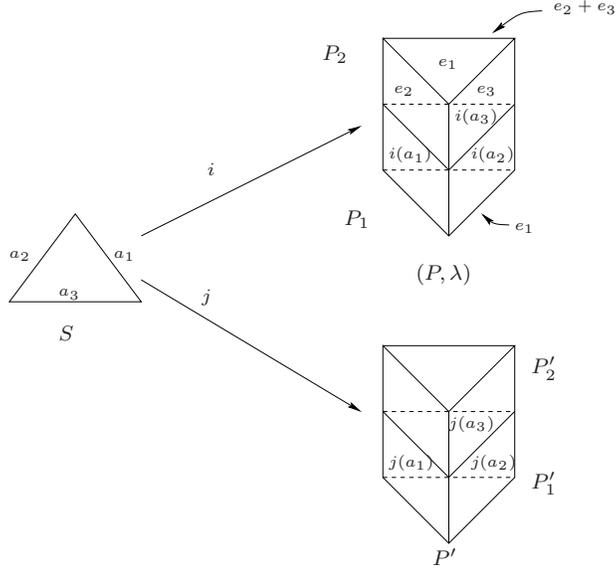
\[   \begin{picture}(0,0)%
\includegraphics{4.pstex}%
\end{picture}%
\setlength{\unitlength}{1816sp}%
\begingroup\makeatletter\ifx\SetFigFont\undefined%
\gdef\SetFigFont#1#2#3#4#5{%
  \reset@font\fontsize{#1}{#2pt}%
  \fontfamily{#3}\fontseries{#4}\fontshape{#5}%
  \selectfont}%
\fi\endgroup%
\begin{picture}(8283,7803)(1189,-8533)
\put(8326,-7411){\makebox(0,0)[lb]{\smash{{\SetFigFont{8}{9.6}{\rmdefault}{\mddefault}{\updefault}$P'_1$}}}}
\put(1876,-5386){\makebox(0,0)[lb]{\smash{{\SetFigFont{8}{9.6}{\rmdefault}{\mddefault}{\updefault}$S$}}}}
\put(7051,-1636){\makebox(0,0)[lb]{\smash{{\SetFigFont{6}{7.2}{\rmdefault}{\mddefault}{\updefault}$e_1$}}}}
\put(8101,-3886){\makebox(0,0)[lb]{\smash{{\SetFigFont{6}{7.2}{\rmdefault}{\mddefault}{\updefault}$e_1$}}}}
\put(7576,-2011){\makebox(0,0)[lb]{\smash{{\SetFigFont{6}{7.2}{\rmdefault}{\mddefault}{\updefault}$e_3$}}}}
\put(3901,-3136){\makebox(0,0)[lb]{\smash{{\SetFigFont{7}{8.4}{\rmdefault}{\mddefault}{\updefault}$i$}}}}
\put(3826,-4861){\makebox(0,0)[lb]{\smash{{\SetFigFont{7}{8.4}{\rmdefault}{\mddefault}{\updefault}$j$}}}}
\put(8626,-886){\makebox(0,0)[lb]{\smash{{\SetFigFont{6}{7.2}{\rmdefault}{\mddefault}{\updefault}$e_2+e_3$}}}}
\put(1876,-4786){\makebox(0,0)[lb]{\smash{{\SetFigFont{6}{7.2}{\rmdefault}{\mddefault}{\updefault}$a_3$}}}}
\put(7276,-2386){\makebox(0,0)[lb]{\smash{{\SetFigFont{6}{7.2}{\rmdefault}{\mddefault}{\updefault}$i(a_3)$}}}}
\put(7201,-6586){\makebox(0,0)[lb]{\smash{{\SetFigFont{6}{7.2}{\rmdefault}{\mddefault}{\updefault}$j(a_3)$}}}}
\put(6451,-2011){\makebox(0,0)[lb]{\smash{{\SetFigFont{6}{7.2}{\rmdefault}{\mddefault}{\updefault}$e_2$}}}}
\put(2626,-4261){\makebox(0,0)[lb]{\smash{{\SetFigFont{6}{7.2}{\rmdefault}{\mddefault}{\updefault}$a_1$}}}}
\put(6376,-2911){\makebox(0,0)[lb]{\smash{{\SetFigFont{6}{7.2}{\rmdefault}{\mddefault}{\updefault}$i(a_1)$}}}}
\put(7501,-2911){\makebox(0,0)[lb]{\smash{{\SetFigFont{6}{7.2}{\rmdefault}{\mddefault}{\updefault}$i(a_2)$}}}}
\put(6376,-7111){\makebox(0,0)[lb]{\smash{{\SetFigFont{6}{7.2}{\rmdefault}{\mddefault}{\updefault}$j(a_1)$}}}}
\put(7501,-7111){\makebox(0,0)[lb]{\smash{{\SetFigFont{6}{7.2}{\rmdefault}{\mddefault}{\updefault}$j(a_2)$}}}}
\put(1201,-4261){\makebox(0,0)[lb]{\smash{{\SetFigFont{6}{7.2}{\rmdefault}{\mddefault}{\updefault}$a_2$}}}}
\put(6751,-4561){\makebox(0,0)[lb]{\smash{{\SetFigFont{8}{9.6}{\rmdefault}{\mddefault}{\updefault}$(P, \lambda)$}}}}
\put(5776,-3811){\makebox(0,0)[lb]{\smash{{\SetFigFont{8}{9.6}{\rmdefault}{\mddefault}{\updefault}$P_1$}}}}
\put(5476,-1561){\makebox(0,0)[lb]{\smash{{\SetFigFont{8}{9.6}{\rmdefault}{\mddefault}{\updefault}$P_2$}}}}
\put(6976,-8461){\makebox(0,0)[lb]{\smash{{\SetFigFont{8}{9.6}{\rmdefault}{\mddefault}{\updefault}$P'$}}}}
\put(8326,-5836){\makebox(0,0)[lb]{\smash{{\SetFigFont{8}{9.6}{\rmdefault}{\mddefault}{\updefault}$P'_2$}}}}
\end{picture}%
\centering
   \]
   So we may obtain a coloring $\bar{\lambda}$ on $S$ by mapping $a_1\longmapsto e_2, a_2\longmapsto e_3, a_3\longmapsto e_2+e_3$. Choose an isomorphism pair $(\psi, \rho)$ as follows: $\psi: \mathcal{F}(S)\longrightarrow\mathcal{F}(S)$ is defined by $a_1\longmapsto a_2, a_2\longmapsto a_3, \mathrm{}a_3\longmapsto a_1$, and $\rho=\begin{bmatrix}
1 & 0 & 0\\
1 & 0 & 1\\
1 & 1 & 1
\end{bmatrix}$. In a similar way to Example~\ref{special example}, we can give two combinatorial equivalences
$g_r: P_r\longrightarrow P'_r$ $(r=1,2)$ and
 a coloring $\lambda'$ on $P'$ as shown in the following picture:
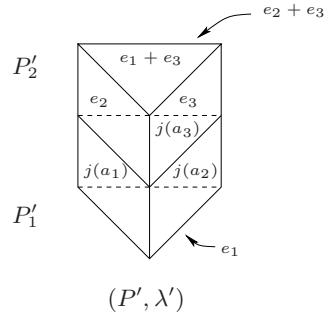
\[   \begin{picture}(0,0)%
\includegraphics{5.pstex}%
\end{picture}%
\setlength{\unitlength}{1973sp}%
\begingroup\makeatletter\ifx\SetFigFont\undefined%
\gdef\SetFigFont#1#2#3#4#5{%
  \reset@font\fontsize{#1}{#2pt}%
  \fontfamily{#3}\fontseries{#4}\fontshape{#5}%
  \selectfont}%
\fi\endgroup%
\begin{picture}(3996,3912)(5476,-4642)
\put(6826,-1486){\makebox(0,0)[lb]{\smash{{\SetFigFont{6}{7.2}{\rmdefault}{\mddefault}{\updefault}$e_1+e_3$}}}}
\put(8101,-3886){\makebox(0,0)[lb]{\smash{{\SetFigFont{6}{7.2}{\rmdefault}{\mddefault}{\updefault}$e_1$}}}}
\put(7576,-2011){\makebox(0,0)[lb]{\smash{{\SetFigFont{6}{7.2}{\rmdefault}{\mddefault}{\updefault}$e_3$}}}}
\put(8626,-886){\makebox(0,0)[lb]{\smash{{\SetFigFont{6}{7.2}{\rmdefault}{\mddefault}{\updefault}$e_2+e_3$}}}}
\put(6451,-2011){\makebox(0,0)[lb]{\smash{{\SetFigFont{6}{7.2}{\rmdefault}{\mddefault}{\updefault}$e_2$}}}}
\put(7276,-2386){\makebox(0,0)[lb]{\smash{{\SetFigFont{6}{7.2}{\rmdefault}{\mddefault}{\updefault}$j(a_3)$}}}}
\put(7501,-2911){\makebox(0,0)[lb]{\smash{{\SetFigFont{6}{7.2}{\rmdefault}{\mddefault}{\updefault}$j(a_2)$}}}}
\put(6376,-2911){\makebox(0,0)[lb]{\smash{{\SetFigFont{6}{7.2}{\rmdefault}{\mddefault}{\updefault}$j(a_1)$}}}}
\put(5476,-3511){\makebox(0,0)[lb]{\smash{{\SetFigFont{9}{10.8}{\rmdefault}{\mddefault}{\updefault}$P'_1$}}}}
\put(5476,-1636){\makebox(0,0)[lb]{\smash{{\SetFigFont{9}{10.8}{\rmdefault}{\mddefault}{\updefault}$P'_2$}}}}
\put(6676,-4561){\makebox(0,0)[lb]{\smash{{\SetFigFont{9}{10.8}{\rmdefault}{\mddefault}{\updefault}$(P', \lambda')$}}}}
\end{picture}%
\centering
   \]
   In this case, we know from \cite[Lemma 5.1]{ly} that the sector manifold $\mathcal{S}$ is a disjoint union of two copies of ${\Bbb R}P^2$. Furthermore, we may conclude from \cite[Corollary 5.2]{ly} that $\widetilde{M}(\lambda)$ is still a small cover over $(P, \lambda)$, and $M(\lambda')$ is a small cover over $(P',\lambda')$. Thus, $\widetilde{M}(\lambda)$ and $M(\lambda')$ are homeomorphic but not (weakly) equivariantly homeomorphic by \cite[Lemma 4.3]{ly}.
\end{example}

\begin{rem}
Clearly the definition of
$\Psi$ depends on  $(\psi,\rho)$ and
$v_0$, so we may write $\Psi$ as $\Psi(\psi,\rho, v_0)$ to indicate
this dependence.
\end{rem}

By Lemma~\ref{induced}, we have known that  $\widetilde{M}(\lambda)$ is identified with $M(\lambda')$ in the sense of homeomorphism.  Now we can use  \cite[Chapter 8, Theorem 2.2]{h} to further analyze when $M(\lambda)$ and $M(\lambda')$ are homeomorphic.

\begin{defn}\label{good twist}
We say that $\Psi(\psi,\rho, v_0)$ is {\em extendable} to  $M_2$ if  there is a  self-homeomorphism
$\widetilde{\Psi}:M_2\rightarrow M_2$ such that
$\widetilde{\Psi}\circ\iota_2=\iota_2\circ\Psi$.
\end{defn}

\begin{thm}\label{extend}
If $\Psi(\psi,\rho, v_0)$ is extendable to
$M_2$, then $M(\lambda)$ is homeomorphic to $M(\lambda')$.
\end{thm}

\begin{proof}
Let $\widetilde{\Psi}:M_2\rightarrow M_2$ be a self-homeomorphism
with $\widetilde{\Psi}\circ\iota_2=\iota_2\circ\Psi$. By
Lemma~\ref{induced}, one may identify $M(\lambda')$ with
$\widetilde{M}(\lambda)=M_1\coprod
M_2/\iota_1(x)\sim\iota_2(\Psi(x))$. Define $H:
M(\lambda)\longrightarrow \widetilde{M}(\lambda)$ by
$$H(z)=\begin{cases} z &
\text{if $z\in M_1\subset M(\lambda)$}\\
\widetilde{\Psi}(z) & \text{if $z\in M_2\subset M(\lambda)$.}
\end{cases}$$
Now, in order to  shows that $H$ is a well-defined
homeomorphism, by \cite[Chapter 8, Theorem 2.2]{h}, it suffices to prove that for $x\in \mathcal{S}$, if
$\iota_1(x)\sim\iota_2(x)$ in $M(\lambda)$, then $H(\iota_1(x))\sim H(\iota_2(x))$ in $\widetilde{M}(\lambda)$.
This is obvious since $H(\iota_1(x))=\iota_1(x)$ and $H(\iota_2(x))=\widetilde{\Psi}(\iota_2(x))=\iota_2(\Psi(x))$.
\end{proof}

\begin{rem}
Take an automorphism $\widetilde{\psi}$ of $P_2$ such that
$\widetilde{\psi}(i_2(S))=i_2(S)$. Then the restriction
$\widetilde{\psi}|_{i_2(S)}$ gives an automorphism of $i_2(S)$. If
we can choose $(\psi, \rho)$ with the property that $\psi\circ
i_2=\widetilde{\psi}|_{i_2(S)}$ and $\rho\circ
\lambda=\lambda\circ\widetilde{\psi}$, then it is easy to check that
$(\widetilde{\psi}, \rho)$ can induce a self-homeomorphism
$\widetilde{\Psi}(\widetilde{\psi}, \rho, v_0)$ of $M_2$, which is
defined by $\{(p, v)\}\longmapsto \{(\widetilde{\psi}(p),
\rho(v)+v_0)\}$. Then $\widetilde{\Psi}(\widetilde{\psi}, \rho,
v_0)$ is an extension of $\Psi(\psi,\rho, v_0)$. In this case, we
see that $M(\lambda)$ is homeomorphic to $M(\lambda')$. However, the
condition $\rho\circ \lambda=\lambda\circ\widetilde{\psi}$ results
in a little bit difficulty for the choice of $(\widetilde{\psi},
\rho)$.
\end{rem}

\cite[Chapter 8, Theorem 2.3]{h} provides us much
insight into a further analysis of $\Psi(\psi,\rho, v_0)$ for
the extension of $\Psi(\psi,\rho, v_0)$.  We say that   $\Psi(\psi,\rho,v_0)$ is called a \emph{good twist} if it is isotopic to the identity.

\begin{thm}\label{t:GENERAL}
If $\Psi(\psi,\rho, v_0)$ is a good twist, then
$M(\lambda)$ is homeomorphic to $M(\lambda')$.
\end{thm}
\begin{proof} Assume that $\Psi(\psi,\rho, v_0)$ is a good twist. Then there is an isotopy $\bar{\Psi}: \mathcal{S}\times[0,1]\rightarrow \mathcal{S}$  such that $\bar{\Psi}(\cdot,0)=\Psi$ and
$\bar{\Psi}(\cdot,1)=\text{id}$.
Since the image of the embedding $i:S\rightarrow P$ does not contain
any vertex of $P$, we can extend $i$ to an embedding $\widetilde{i}:
S\times[-1,1]\rightarrow P$ such that $i=\widetilde{i}(\cdot,0)$ and
each $\widetilde{i}(\cdot, t)$ is a sector. We can further assume
that $\widetilde{i}(S\times[-1,0])\subset P_1$ and
$\widetilde{i}(S\times[0,1])\subset P_2$. In the world of
topology, $\widetilde{i}$ corresponds to an embedding
$\widetilde{\iota}: \mathcal{S}\times[-1,1]\rightarrow M(\lambda)$
such that $\iota=\widetilde{\iota}(\cdot,0)$, and
$\widetilde{\iota}(\mathcal{S}\times[0,1])\subset M_2$.  Now define
$\widetilde{\Psi}:M_2\rightarrow M_2$ by
$$\widetilde{\Psi}(y)=\begin{cases} \widetilde{\iota}(\bar{\Psi}(x,t),t) &
\text{if
$y=\widetilde{\iota}(x,t)\in \text{Im}\widetilde{\iota}$}\\
y & \text{if $y\not\in \text{Im}\widetilde{\iota}$.}
\end{cases}$$
We check easily that $\widetilde{\Psi}$ is a self-homeomorphism of
$M_2$ such that
 $\widetilde{\Psi}|_{\iota_2(\mathcal{S})}=\iota_2\circ\Psi$. Moreover, Theorem~\ref{t:GENERAL} follows by
 applying Theorem~\ref{extend}.
\end{proof}

 We see from the proof of Lemma~\ref{induced} that the homeomorphism type of
$M(\lambda')$ doesn't depend on the choice of $v_0$. So to apply
Theorem~\ref{t:GENERAL} we can choose a suitable $v_0$ such that
$\Psi(\psi,\rho,v_0)$ meets the required condition.

\begin{rem}
The sector method above provides a way of how to construct a
homeomorphism between two small covers $M(\lambda)\longrightarrow P$
and $M(\lambda')\longrightarrow P'$ regardless of whether $P$ is
combinatorially equivalent to $P'$ or not. In addition, the sector
method also gives an approach of how to construct a new colored
polytope $(P', \lambda')$ from $(P, \lambda)$ by the
isomorphism pair $(\psi, \rho)$ of $S$.
\end{rem}

\section{Application to prisms: Rectangular Sector Method}
\label{Rectangular}

The objective of this section is to give the application of the
sector method to prisms.

Let $\Pol^3(m)$ denote a 3-dimensional prism that is the product of
$[0,1]$ and an $m$-gon  where $m\geq 3$. When $m\neq 4$, let $c,f$
(the \emph{ceiling} and the \emph{floor}) be the two 2-faces of
$\Pol^3(m)$ that are $m$-gons. For the 3-cube (i.e., $m=4$), we
specify two opposite 2-faces and distinguish them as ceiling and
floor. For convenience,  we identify  other 2-faces (i.e., side
2-faces) with $s_1,..., s_m$ in the natural way. Let
$$\Lambda(\Pol^3(m))=\{\lambda\big| \lambda \text{ is a $(\Z_2)^3$-coloring on $\Pol^3(m)$}\}.$$

\subsection{Rectangular sector}
 Generally, any polygon
can become a sector in the setting of all 3-polytopes.  However,
here we shall pay attention on rectangular sectors because this will
be sufficient enough to the classification of all small covers over
prisms.

 Throughout the following, choose  the rectangle
$S=\{(x,y)\in \R^2\big| |x|\leq 1, |y|\leq 1\}$ in a plane $\R^2$.
Clearly $S$ can always be embedded as a sector in any
$(\Z_2)^3$-colored simple 3-polytope $(P^3, \lambda)$. Fix
$\{e_1,e_2,e_3\}$ as a basis of $(\Z_2)^3$, then it is easy to see
that up to Davis-Januszkiewicz equivalence, all possible derived
colorings $\bar{\lambda}:\mathcal{F}(S)\longrightarrow (\Z_2)^3$ and
corresponding sector manifolds $\mathcal{S}$ can be stated as
follows:
\begin{center}
\begin{tabular}{|c|c|c|c|c|c|c|}
\hline  & top edge  & left edge  & bottom edge  & right edge & sector manifold $\mathcal{S}$\\
\hline
   $\bar{\lambda}_1$   & $e_1$ & $e_2$ & $e_1$ & $e_2 $ &  union of 2 tori    \\
    $\bar{\lambda}_2$ & $e_1$ & $e_2$ & $e_1+e_2$ & $e_2$ &  union of 2 Klein bottles    \\
  $\bar{\lambda}_3$     & $e_1$ & $e_2$ & $e_3$ & $e_2$  & torus      \\
    $\bar{\lambda}_4$        & $e_3$ & $e_1$ & $e_1+e_3$ &$e_2$  & Klein bottle \\
  $\bar{\lambda}_5$    & $e_3$ & $e_1$ & $e_1+e_2+e_3$ & $e_2$ & torus \\
\hline
\end{tabular}
\end{center}

 It is well-known that the symmetric group $
\mathrm{Aut}(S)$ of $S$ as a 4-gon is isomorphic to the dihedral
group $\mathcal{D}_4$ of order 8, which  just contains four
reflections. Clearly, each reflection of $S$ may be expressed as a
matrix. For example, the reflection along $y$-axis can be written as
$\text{diag}(-1, 1)$, and the reflection along $x$-axis can be
written as $\text{diag}(1, -1)$.

\subsection{Construction of new colored polytopes from  $(\Pol^3(m),
\lambda)$} Given a pair $(\Pol^3(m), \lambda)$ in
$\Lambda(\Pol^3(m))$. We use the convention that all embedded
rectangular sectors of $(\Pol^3(m), \lambda)$ used here are always
orthogonal to the ceiling and floor of $\Pol^3(m)$.
 Given such a sector $i: S\longrightarrow
\Pol^3(m)$ (note: here we don't need that $i$ must map the top edge
and the bottom edge of $S$ into the ceiling and the floor of
$\Pol^3(m)$, respectively), it is easy to see that there are two
side faces $s_k, s_l$ $(k\neq l$ and $k<l)$ such that
$S\longrightarrow \Pol^3(m)$ is essentially determined by $s_k$ and
$s_l$ and it is called the sector \emph{at} $s_k, s_l$,
denoted by $i(k,l)$. Let $P_1$ and $P_2$ be two prisms cut out by $i(k,l):
S\longrightarrow \Pol^3(m)$ from $\Pol^3(m)$. Without the loss of generality, assume that $P_2$ contains side
faces $s_k,...,s_l$ of $\Pol^3(m)$.

Now, let us use the sector method to discuss how to construct  new
colored 3-polytopes $(\Pol',\lambda')$ from $(\Pol^3(m),\lambda)$ by
 isomorphism pairs $(\psi,\rho)$ at the sector $i(k,l):
S\longrightarrow \Pol^3(m)$. For our purpose, we wish that (1)
$\Pol'$ is still combinatorially equivalent to $\Pol^3(m)$, and (2)
 $\Psi(\psi,\rho, v_0)$ is a good twist, so that each $M(\lambda')$
is homeomorphic to $M(\lambda)$ by Theorem~\ref{t:GENERAL}. This
will depend upon the choices of $\psi$ and $\rho$. Actually, the
construction of $\Pol'$ depends upon the choice of $\psi$, and the
definition of $\lambda'$ depends upon the choice of $\rho$. In
particular, $v_0$ will provide a convenience for the choices of
$(\psi, \rho)$.

\vskip .2cm {\em Convention}: $\psi=\text{\rm diag}(-1,1)$ or
$\text{\rm diag}(1,1)$ means that $i(k,l)$  maps the top edge and
the bottom edge of $S$ into the ceiling and the floor of
$\Pol^3(m)$, respectively, and $\psi=\text{\rm diag}(1,-1)$ means
that $i(k,l)$ maps the left edge and the right edge of $S$ into the
ceiling and the floor of $\Pol^3(m)$, respectively.

\begin{lem}\label{construct}
If $\psi=\text{\rm diag}(-1,1)$ or $\text{\rm diag}(1,-1)$ or
  $\text{\rm diag}(1,1)$, then we can construct  a
new polytope $\Pol'$  from $\Pol^3(m)$ by
 $\psi$ at the sector $i(k,l):
S\longrightarrow \Pol^3(m)$ such that  $\Pol'$ is combinatorially
equivalent to $\Pol^3(m)$.
\end{lem}
\begin{proof} The case $\psi=\text{\rm diag}(1,1)$ is trivial.
Actually, in this case we can just choose $P'_r=P_r$ and let $f_r:
P_r\longrightarrow P'_r$ be the identity, where $r=1,2$.  If
$\psi=\text{diag}(-1,1)$, to construct $\Pol'$, we first choose
$P'_1=P_1$  and $f_1$ as the  identity from $P_1\longrightarrow
P'_1$, and then choose $P'_2$ as the image of mirror reflection $R$
of $P_2$ along a 2-plane $H$ orthogonal to the ceiling and floor of
$P_2$ with $H\cap P_2=\emptyset$ (i.e., intuitively $P'_2$ is
obtained by reserving the ordering of the side faces $s_k, ..., s_l$
of $P_2$) and $f_2$ as the  homeomorphism induced by the reflection
$R$, as shown in the following figure.
 
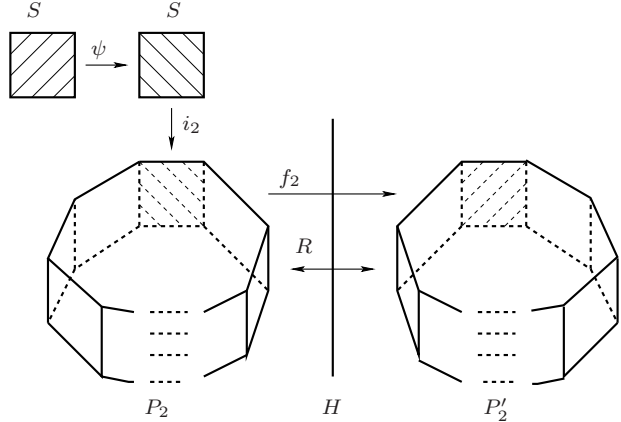
\[   \begin{picture}(0,0)%
\includegraphics{p.pstex}%
\end{picture}%
\setlength{\unitlength}{1776sp}%
\begingroup\makeatletter\ifx\SetFigFont\undefined%
\gdef\SetFigFont#1#2#3#4#5{%
  \reset@font\fontsize{#1}{#2pt}%
  \fontfamily{#3}\fontseries{#4}\fontshape{#5}%
  \selectfont}%
\fi\endgroup%
\begin{picture}(8541,5823)(868,-6367)
\put(4876,-4036){\makebox(0,0)[lb]{\smash{{\SetFigFont{8}{9.6}{\rmdefault}{\mddefault}{\updefault}$R$}}}}
\put(3301,-2311){\makebox(0,0)[lb]{\smash{{\SetFigFont{8}{9.6}{\rmdefault}{\mddefault}{\updefault}$i_2$}}}}
\put(4651,-3061){\makebox(0,0)[lb]{\smash{{\SetFigFont{8}{9.6}{\rmdefault}{\mddefault}{\updefault}$f_2$}}}}
\put(2026,-1261){\makebox(0,0)[lb]{\smash{{\SetFigFont{8}{9.6}{\rmdefault}{\mddefault}{\updefault}$\psi$}}}}
\put(3076,-736){\makebox(0,0)[lb]{\smash{{\SetFigFont{8}{9.6}{\rmdefault}{\mddefault}{\updefault}$S$}}}}
\put(1126,-736){\makebox(0,0)[lb]{\smash{{\SetFigFont{8}{9.6}{\rmdefault}{\mddefault}{\updefault}$S$}}}}
\put(2776,-6286){\makebox(0,0)[lb]{\smash{{\SetFigFont{8}{9.6}{\rmdefault}{\mddefault}{\updefault}$P_2$}}}}
\put(5251,-6286){\makebox(0,0)[lb]{\smash{{\SetFigFont{8}{9.6}{\rmdefault}{\mddefault}{\updefault}$H$}}}}
\put(7501,-6286){\makebox(0,0)[lb]{\smash{{\SetFigFont{8}{9.6}{\rmdefault}{\mddefault}{\updefault}$P'_2$}}}}
\end{picture}%
\centering
   \]
Now, we clearly see that $\Pol'$ can be defined as $P_1\coprod P'_2/
i_1(s)\sim f_2\circ i_2\circ \psi(s)$, which is combinatorially
equivalent to $\Pol^3(m)$. In a similar way, we can prove the case
$\psi=\text{diag}(1,-1)$.
\end{proof}

\subsection{Good twists} Now suppose that $\Pol'$, which is just constructed from
$\Pol^3(m)$ by a  $\psi$ at the sector $i(k,l): S\longrightarrow
\Pol^3(m)$,  is combinatorially equivalent to $\Pol^3(m)$. Then, as
stated in Section~\ref{sector},  we can use $\rho$ to give a
coloring $\lambda'$ on $\Pol'$ as long as $\rho$ satisfies the
equation $\rho\circ \bar{\lambda}=\bar{\lambda}\circ \psi$. To
guarantee that $M(\lambda')$ is homeomorphic to $M(\lambda)$, we
need choose $(\psi, \rho)$ carefully such that $\Psi(\psi,\rho,
v_0)$ is a good twist.

 Based upon the possible values of $\bar{\lambda}$ (see
the table above),  we may find some good twists and list them as
follows:
\begin{center}
\begin{tabular}{|c|c|c|c|ccc|c|c|}
\hline Sector &$(S, \bar{\lambda})$ & $\Psi(\psi,\rho,v_0)$ & $\psi$ &$\rho(e_1)$ &$\rho(e_2)$ & $\rho(e_3)$ & $v_0$\\
\hline
    $S(1)$    & $(S, \bar{\lambda}_1)$    &$\Psi(\psi,\rho,v_0)$         & $\text{diag}(1,-1)$ & $e_1$ & $e_2$ & $e_3$ & $e_1$ \\
    $S(2_1)$  & $(S, \bar{\lambda}_2)$    &$\Psi(\psi,\rho,v_0)$   & $\text{diag}(1,1) $ & $e_1$ & $e_2$ & $e_3+e_2$ & $0$ \\
    $S(2_2)$  & $(S, \bar{\lambda}_2)$    &$\Psi(\psi,\rho,v_0)$   & $\text{diag}(1,-1)$ & $e_1+e_2$& $e_2$ & $e_3 $ & $e_1$  \\
    $S(3_1)$  & $(S, \bar{\lambda}_3)$        &$\Psi(\psi,\rho,v_0)$             & $\text{diag}(-1,1)$ & $e_1$ & $e_2$ & $e_3$ & $e_2$\\
    $S(3_2)$  & $(S, \bar{\lambda}_3)$        &$\Psi(\psi,\rho,v_0)$            & $\text{diag}(1,-1)$ & $e_3$ & $e_2$ & $e_1$ & $e_1$\\
    $S(4)$    & $(S, \bar{\lambda}_4)$   &$\Psi(\psi,\rho,v_0)$      & $\text{diag}(1,-1)$ & $e_1$ & $e_2$ & $e_3+e_1$& $e_3$\\
    $S(5)$    & $(S, \bar{\lambda}_5)$&$\Psi(\psi,\rho,v_0)$            & $\text{diag}(-1,1)$ & $e_2$ & $e_1$ & $e_3$& $e_1$\\
\hline\end{tabular}
\end{center}

\begin{rem}
It should be pointed out that we have not listed all such good
twists.  At least we omit compositions of good twists that have
already appeared in the above  table. However, as we shall see,
those good twists listed above are sufficient in the further
applications.
\end{rem}

Here we only give a detailed argument of  $S(1)$  because all other
cases can be checked  similarly. On $S(1)$, we know that
$\mathcal{S}={[-1,1]}^2\times(\Z_2)^3/\sim$ is the disjoint union of
two  tori. Then  we can write $\mathcal{S}=\{(z_1,z_2,\alpha)\big|
z_k\in S^1\subset{\Bbb C}, k=1, 2,  \alpha\in\{0,1\}\}$, such that
for
$\mathbf{x}=\{((x_1,x_2),v=a_1e_2+a_2e_1+a_3e_3)\}\in\mathcal{S}$,
there is the following one-one correspondence
$$z_k(\mathbf{x})=\begin{cases}\mathrm{exp}(\mathbf{i}x_k\pi/2) &\text{
if }a_k=0\\
\mathrm{exp}(\mathbf{i}(\pi-x_k\pi/2)) &\text{ if
}a_k\neq0\end{cases}$$ for $k=1,2$ and $\alpha(\mathbf{x})=a_3$. Now
we consider the map $\Psi(\psi,\rho,v_0)$ where
$\psi(x_1,x_2)=(x_1,-x_2)$, $\rho=\text{id}$ and $v_0=e_1$. An easy
computation yields that $\Psi(z_1, z_2, \alpha)=(z_1, -z_2,
\alpha)$, which is clearly isotopic to the identity via the homotopy
$((z_1$, $z_2$, $\alpha), t)$$ \mapsto $$(z_1$,
$z_2\mathrm{exp}({\mathbf{i}\pi t})$, $\alpha)$ where $t\in[0,1]$.
Thus $\Psi(\psi,\rho,v_0)$ is a good twist.

\section{Operations on coloring sequences and canonical forms}\label{prism}

Now we apply the developed rectangular sector method to study  small
covers over prisms.

Given a pair $(\Pol^3(m), \lambda)$ in $\Lambda(\Pol^3(m))$ and a
sector $i(k,l): S\longrightarrow \Pol^3(m)$. We have known how to
construct a $(\Pol',\lambda')$ from $(\Pol^3(m),\lambda)$ by an isomorphism pair $(\psi,\rho)$ at the sector $S\longrightarrow
\Pol^3(m)$.
 Indeed, by Lemma~\ref{construct}, if $\psi=\text{id}$ or $\text{\rm
diag}(-1,1)$ or $\text{\rm diag}(1,-1)$, then  $\Pol'$ is also a
$\Pol^3(m)$ with the same ceiling and floor coloring, and $\lambda'$
has the same  side coloring sequence  as $\lambda$ on sides faces
from $s_{l+1}$ to $s_{k-1}$. For $k\leq r\leq l$,  if
$\psi=\text{id}$, then $\lambda'(s_r)=\rho^{-1}\lambda(s_r)$; if
$\psi=\text{\rm diag}(-1,1)$ or $\text{\rm diag}(1,-1)$, then
$\lambda'(s_r)=\rho^{-1}\lambda(s_{k+l-r})$, that is, we reflect the
sequence from $s_k$ to $s_l$ and apply the linear transformation
$\rho$. By Theorem~\ref{t:GENERAL}, when the derived coloring
$\bar{\lambda}$ of the sector at $s_k,s_l$ and $(\psi,\rho)$ match a
case in the table of last section, we can conclude that
$M(\lambda')$ is homeomorphic to $M(\lambda)$.
  Thus we can reduce $(\Pol^3(m), \lambda)$ to
$(\Pol^3(m), \lambda')$ without changing the homeomorphism type of
the small cover. In this case,  both $(\Pol^3(m), \lambda)$ and
$(\Pol^3(m), \lambda')$ are said to be {\em sector-equivalent},
denoted by $(\Pol^3(m), \lambda)\approx (\Pol^3(m), \lambda')$ or
simply $\lambda\approx\lambda'$.

Based upon this rectangular sector method, we shall show that
$\Lambda(\Pol^3(m))$ contains some basic colored pairs, called
``canonical forms'', such that any pair in $\Lambda(\Pol^3(m))$ is
sector-equivalent  to one of canonical forms. This means that up to
homeomorphism,  those canonical forms determine all small covers
over $\Pol^3(m)$.

For a convenience, after fixing the colorings of ceiling and floor,
we use the {\em convention} that a coloring on $\Pol^3(m)$ will
simply be described as a sequence by writing its side face colorings
in order, keeping in mind that the first one is next to the last.

\begin{defn}A coloring
$\lambda\in\lambda(\Pol^3(m))$ is said to be \emph{2-independent} if
all $\lambda(s_i), i=1,..., m$, span a 2-dimensional subspace of
$(\Z_2)^3$; otherwise it's said to be \emph{3-independent}. If
$\lambda(c)=\lambda(f)$, then $\lambda$ is said to be
\emph{trivial}; otherwise \emph{nontrivial}. \end{defn}

The argument is divided into two cases: (i) $\lambda$ is trivial;
(ii)  $\lambda$ is nontrivial.

\subsection{Trivial colorings} Given a pair $(\Pol^3(m), \lambda)$ in
$\Lambda(\Pol^3(m))$, throughout the following suppose that
$\lambda$ is trivial with $\lambda(c)=\lambda(f)=e_1$. Let
$\{\alpha,\beta,e_1\}$ be a basis of $(\Z_2)^3$, and let
$\gamma=\alpha+\beta$. Write $\bar{\alpha}=\alpha+e_1$,
$\bar{\beta}=\beta+e_1$ and $\bar{\gamma}=\gamma+e_1$. We say that
$\lambda$ satisfies the {\em property $(\star)$} if all three letters
$\alpha, \beta, \gamma$ (with or without bar we don't care) appear
in its coloring sequence.

\vskip .2cm Applying sectors $S(1), S(2_1), S(2_2)$ and $S(3_2)$ to
the trivial coloring $\lambda$ gives the following four fundamental
operations on its coloring sequence $(\lambda(s_1), ...,
\lambda(s_m))$:

\begin{enumerate}
\item[$\text{O}_1$] Take  two side faces $s_k,s_l$ $(k<l)$ with the same coloring and then  use $S(1)$
to reflect the coloring sequence of $s_k, s_{k+1}, ...,  s_l$.
\item[$\text{O}_{21}$] Take two faces $s_k, s_l$ $(k<l)$  with $\lambda(s_k)=\lambda(s_l)+\lambda(c)$
(without loss of generality,  assume that
$\{\lambda(s_k),\lambda(s_l)\}=\{\alpha,\bar{\alpha}\}$), and then
by using $S(2_1)$, we can do a linear transform
$(e_1,\alpha,\beta)\mapsto(e_1,\alpha,\bar{\beta})$ to change the
coloring sequence of $s_k, s_{k+1},  ...,  s_l$.
\item[$\text{O}_{22}$]
Take two faces $s_k, s_l$ $(k<l)$  with
$\lambda(s_k)=\lambda(s_l)+\lambda(c)$ and
$\{\lambda(s_k),\lambda(s_l)\}=\{\alpha,\bar{\alpha}\}$ as above,
then by using $S(2_2)$ we can reflect the coloring sequence of $s_k,
s_{k+1},  ...,  s_l$ and do a linear transform
$(e_1,\alpha,\beta)\mapsto(e_1,\bar{\alpha},\beta)$ to change the
reflected coloring sequence.
\item[$\text{O}_{32}$] Take $s_k,s_l$ $(k<l)$  with $\lambda(s_k),\lambda(s_l),e_1$
independent, and then by using $S(3_2)$ we can reflect the coloring
sequence of $s_k, s_{k+1},.., s_l$ and   do a linear transform
$(\lambda(s_k),\lambda(s_l),e_1)\mapsto(\lambda(s_l),\lambda(s_k),e_1)$
to change the reflected coloring sequence.
\end{enumerate}

\begin{lem}\label{trivial}
The trivial coloring $\lambda$ with  the property $(\star)$ is
always sector-equivalent to a coloring whose coloring sequence
contains only one of both $\mathrm{\gamma}$ and
$\mathrm{\bar{\gamma}}$.
\end{lem}

\begin{proof} Let $\widetilde{\gamma}=\gamma$ or $\bar{\gamma}$.
With no loss, assume that the time number $\ell$ of
$\widetilde{\gamma}$ appearing  in the coloring sequence of
$\lambda$ is greater than one. Up to Davis-Januszkiewicz
equivalence, one also may assume that $\ell< m/2$. By the definition
of $\lambda$, it is easy to see that  any two $\widetilde{\gamma}$'s
in the coloring sequence  cannot become neighbors. Let
$\widetilde{\gamma}, x_1,...,x_r, \widetilde{\gamma}, y$ with $x_i,
y\not=\widetilde{\gamma}$ be a subsequence of the coloring sequence.
If $r>1$, we proceed as follows:

\begin{enumerate}
\item  When $x_1=y$, by doing the operation $\text{O}_1$ on $x_1,..., x_r, \widetilde{\gamma},y$,  we may only change the
subsequence $\widetilde{\gamma}, x_1,..., x_r, \widetilde{\gamma},y$
into  $\widetilde{\gamma}, y, \widetilde{\gamma}, x_r,...,  x_1$ in
the coloring sequence, and the value of $\ell$ is unchanged.
\item When $x_1-y=e_1$, with no loss one may assume that
$x_1={\alpha},y={\bar{\alpha}}$. Then by doing the operation
$\text{O}_{22}$ on $x_1=\alpha, x_2, ..., x_r,
\widetilde{\gamma},y={\bar{\alpha}}$, we may only change the
subsequence $\widetilde{\gamma}, \alpha, x_2,..., x_r,
\widetilde{\gamma},{\bar{\alpha}}$ into $\widetilde{\gamma}, \alpha,
\widetilde{\gamma}, x'_r,...,  x'_2, \bar{\alpha}$ with
$x'_i\not=\widetilde{\gamma}$, and the value of $\ell$ is unchanged.
\item When $x_1,y,e_1$ are linearly independent, with no loss one  may
assume that $x_1=\alpha, y=\beta$. Then by doing the operation
$\text{O}_{32}$ on $x_1=\alpha, x_2, ..., x_r$, $\widetilde{\gamma},
y=\beta$,  we may only change the subsequence $\widetilde{\gamma},
\alpha, x_2,..., x_r, \widetilde{\gamma}, \beta$ into
$\widetilde{\gamma}$, $\alpha$, $\widetilde{\gamma}, x'_r,..., x'_2,
\beta$ with $x'_i\not=\widetilde{\gamma}$, and the value of $\ell$
is unchanged.
\end{enumerate}
Thus, we may reduce the coloring $\lambda$ to another coloring with
the following coloring sequence
\begin{equation}\label{seq}
(\widetilde{\gamma}, y_1, \widetilde{\gamma}, y_2, ...,
\widetilde{\gamma}, y_{\ell-1}, \widetilde{\gamma}, y_\ell, z_1,...,
z_{m-2\ell}) \text{ with } m-2\ell>0.\end{equation}
 Without loss of generality, assume that
$y_{\ell-1}=\alpha$. If $y_\ell=\beta$ or $\bar{\beta}$, by doing
the operation $\text{O}_{32}$ on $\widetilde{\gamma}, y_{\ell-1},
\widetilde{\gamma}, y_\ell$, one may change $\widetilde{\gamma},
y_{\ell-1}, \widetilde{\gamma}, y_\ell$ into $ \widetilde{\gamma},
y_\ell, y_{\ell-1}, y_\ell$, so that the coloring
sequence~(\ref{seq})
 is reduced to $(\widetilde{\gamma}$, $y_1$,
$\widetilde{\gamma}$, $y_2, ...$, $\widetilde{\gamma}, y_{\ell-2}$,
$\widetilde{\gamma}$, $y_\ell, y_{\ell-1}, y_\ell$, $z_1,...,
z_{m-2\ell})$. If $y_\ell=\alpha$ or $\bar{\alpha}$, then
$z_1=\beta$ or $\bar{\beta}$. By doing the operation $\text{O}_{32}$
on $\widetilde{\gamma}, y_{\ell-1}, \widetilde{\gamma}, y_\ell,
z_1$, one may change $\widetilde{\gamma}, y_{\ell-1},
\widetilde{\gamma}, y_\ell, z_1$ into $\widetilde{\gamma}, y_{\ell},
z_1, y_{\ell-1}, z_1$, so that the coloring sequence~(\ref{seq}) is
reduced to $(\widetilde{\gamma}$, $y_1$, $\widetilde{\gamma}$, $y_2,
...$, $\widetilde{\gamma}, y_{\ell-2}$, $\widetilde{\gamma}$,
$y_\ell, z_1, y_{\ell-1}$, $z_1,..., z_{m-2\ell})$. So we have
managed to reduce the number $\ell$ of $\widetilde{\gamma}$'s by 1.
We can continue this process until we reach $\ell=1$, as desired.
\end{proof}

Now let us determine the ``canonical form'' of the trivial coloring
$\lambda$ on $\Pol^3(m)$.

First let us consider the case in which $\lambda$ is 2-independent.

\begin{prop}\label{c1}
Suppose that $\lambda$ is 2-independent. Then
\begin{enumerate}
\item If $\lambda$ doesn't possess the property $(\star)$, then $m$ is even and $\lambda$ is sector-equivalent to
the canonical form $\lambda_{C_1}$ with the coloring sequence
$\mathcal{C}_1=(\alpha$, $ \beta, $$..., \alpha, \beta)$.
\item If $\lambda$  possesses the property $(\star)$,
then $\lambda$ is sector-equivalent to one of  the following two
canonical forms: $\mathrm{(a)}$ $\lambda_{C_2}$ with $m$ even and
with the coloring sequence $\mathcal{C}_2=(\alpha, \gamma, \alpha,
\beta, ..., \alpha, \beta)$; $\mathrm{(b)}$
 $\lambda_{C_3}$ with $m$ odd and with the coloring sequence
$\mathcal{C}_3=(\alpha, \gamma, \beta, \alpha, \beta, ..., \alpha,
\beta)$.
\end{enumerate}
\end{prop}
\begin{proof}
If $\lambda$ doesn't possess the property $(\star)$, then it is easy
to see that $\lambda$ is unique up to Davis-Januszkiewicz
equivalence and $m$ is even. So Proposition~\ref{c1}(1) follows from
this. By Lemma~\ref{trivial}, an easy observation shows that
Proposition~\ref{c1}(2) holds.
\end{proof}

Next let us consider the case in which $\lambda$ is 3-independent.

\begin{prop}\label{c2}
If $\lambda$ is 3-independent without the property $(\star)$, then
$m$ is even and  $\lambda$ is sector-equivalent to  the canonical
form $\lambda_{C_4}$ with the coloring sequence
$\mathcal{C}_4=(\bar{\alpha}, \beta, \alpha, \beta, ..., \alpha,
\beta)$.
\end{prop}

\begin{proof}
Without loss of generality assume that each element in the coloring
sequence $\mathcal{C}$ of $\lambda$ is in the set
$\{\alpha,\beta,\bar{\alpha},\bar{\beta}\}$ and that both $\alpha$
and $\bar{\alpha}$ must appear in $\mathcal{C}$. Since $\alpha$ and
$\bar{\alpha}$ (or $\beta$ and $\bar{\beta}$) can become neighbors,
one has that $m$ must be even.

Similarly to the argument of Lemma~\ref{trivial}, by using the
operations $\text{O}_1$ and $\text{O}_{22}$, we may reduce $\lambda$
to a coloring with the coloring sequence
\begin{equation}\label{seq1}
(\bar{\alpha}, x_1,..., \bar{\alpha}, x_r, \alpha, x_{r+1}, ...,
{\alpha}, x_{{m\over 2}})\end{equation} where $r\geq 1$ and
$x_i=\beta$ or $\bar{\beta}$, and this reduction doesn't change the
number of bars on $\alpha$'s. If $r>1$, by doing the operation
$\text{O}_{22}$ on $\bar{\alpha}, x_{r-1},\bar{\alpha}, x_r,
\alpha$, we may reduce the sequence~(\ref{seq1}) to $$(\bar{\alpha},
x_1,..., \bar{\alpha}, x_{r-2}, \bar{\alpha}, x_r, \alpha, x_{r-1},
\alpha, x_{r+1}, ..., {\alpha}, x_{{m\over 2}}),$$ reducing the
number of bars on $\alpha$'s by one.
 This
process can be carried out  until the sequence~(\ref{seq1}) is
reduced to
\begin{equation}\label{seq2}
(\bar{\alpha}, x_r, \alpha, x_1, ..., {\alpha}, x_{r-2}, \alpha,
x_{r-1}, \alpha, x_{r+1}, ..., {\alpha}, x_{{m\over
2}}).\end{equation} Next, we claim that by using the operation
$\text{O}_{21}$, we may remove all possible bars on $\beta$'s in the
sequence~(\ref{seq2}). In fact, if $x_r=\bar{\beta}$, then applying
the operation $\text{O}_{21}$ on $\bar{\alpha}, x_r, \alpha$, we may
remove the bar on $x_r$. Generally, with no loss, one may assume
that $x_j=\bar{\beta}$ and $x_r=x_l=\beta$ where $l\in\{1,...,
j-1\}$ if $j\leq r+1$ and $j\not=r$,  and $l\in\{1,...,r-1,
r+1,...,j-1\}$ if $j>r+1$. Applying the operation $\text{O}_{21}$ on
$\bar{\alpha}, x_r, ...,  \alpha, x_j=\bar{\beta}, \alpha$, one may
remove the bar on $x_j=\bar{\beta}$, but add the bar on $x_r$ and
$x_l$'s. Again applying the operation $\text{O}_{21}$ on
$\bar{\alpha}, \bar{x}_r, ..., \alpha, \bar{x}_l,..., \alpha,
\bar{x}_{j-1}, \alpha$, one may remove all bars on $\bar{x}_r$ and
$\bar{x}_l$'s. Thus, by carrying on this procedure,  the above claim
holds. This completes the proof.
\end{proof}

\begin{prop}\label{c3}
If $\lambda$ is 3-independent with the property $(\star)$ and $m>4$,
then $\lambda$ is sector-equivalent to one of  the following six
canonical forms:
\begin{enumerate}
\item  $\lambda_{C_5}$ with $m$ odd and with the coloring sequence $\mathcal{C}_5=(\bar{\gamma},\alpha,\beta,...,\alpha,\beta).$
\item  $\lambda_{C_6}$ with $m$ odd and with the coloring sequence $\mathcal{C}_6=(\bar{\gamma},\bar{\alpha},\beta,\alpha,\beta,...,\alpha,\beta).$
\item  $\lambda_{C_7}$ with $m$ odd and with the coloring sequence $\mathcal{C}_7=({\gamma},\bar{\alpha},\beta,\alpha,\beta,...,\alpha,\beta).$
\item  $\lambda_{C_8}$ with $m$ even and with the coloring sequence $\mathcal{C}_8=(\bar{\alpha},\gamma,\bar{\alpha},\beta,\alpha,\beta,...,\alpha,\beta).$
\item  $\lambda_{C_9}$ with $m$ even and with the coloring sequence $\mathcal{C}_9=(\alpha,\gamma,\bar{\alpha},\beta,\alpha,\beta,...,\alpha,\beta).$
\item  $\lambda_{C_{10}}$ with $m$ even and with the coloring sequence $\mathcal{C}_{10}=(\alpha,\bar{\gamma},\alpha,\beta,...,\alpha,\beta).$
\end{enumerate}
\end{prop}

\begin{proof}
By Lemma~\ref{trivial}, one may assume that
$\widetilde{\gamma}(=\gamma$ or $\bar{\gamma})$ appears only one
time in the coloring sequence $\mathcal{C}$ of $\lambda$.

{\bf  Case (I):} {\em $m$ is odd.} If only one of both $\alpha$ and
$\bar{\alpha}$ appears in $\mathcal{C}$ and the same thing also
happens for both $\beta$ and $\bar{\beta}$, then
$\widetilde{\gamma}$ must be $\bar{\gamma}$ so $\mathcal{C}$ is just
$\mathcal{C}_5$.  Otherwise we can carry out our argument as in
Proposition~\ref{c2} on the subsequence in $\mathcal{C}$ of
containing no $\widetilde{\gamma}$, so that $\mathcal{C}$ is
sector-equivalent to $\mathcal{C}_6$ or $\mathcal{C}_7$.

{\bf  Case (II):} {\em $m$ is even.} Consider two neighbors of
$\widetilde{\gamma}$, since $m$ is even, such two neighbors must be
the same letter.  Up to Davis-Januszkiewicz equivalence, one may
assume that they are $\{\alpha, {\alpha}\}$ or
$\{\alpha,\bar{\alpha}\}$.

 If two neighbors of $\widetilde{\gamma}$ are
$\{\alpha,\bar{\alpha}\}$, with no loss assume that
$\mathcal{C}=(\alpha, \widetilde{\gamma}, \bar{\alpha},
x_4,...,x_m)$. Then we can carry out our argument as in
Proposition~\ref{c2} on $ \bar{\alpha}, x_4,...,x_m$, so that
$\mathcal{C}$ may be reduced to
$\mathcal{C}'=(\alpha,\widetilde{\gamma},\bar{\alpha},\beta,\alpha,\beta,...,\alpha,\beta)$.
If $\widetilde{\gamma}=\gamma$, then $\mathcal{C}'$ is just
$\mathcal{C}_8$. If $\widetilde{\gamma}=\bar{\gamma}$, applying the
operation $\text{O}_{22}$ on $\alpha,\bar{\gamma},\bar{\alpha}$, one
may further reduce $\mathcal{C}'$ to $\mathcal{C}_9$.

Now suppose that two neighbors of $\widetilde{\gamma}$ are
$\{\alpha, {\alpha}\}$. If $\mathcal{C}$ only contains $\alpha$ and
$\beta$ except for $\widetilde{\gamma}$, then $\mathcal{C}$ is just
$\mathcal{C}_{10}$. Otherwise, with no loss assume that
$\mathcal{C}$ also contains $\bar{\alpha}$. Up to
Davis-Januszkiewicz equivalence, by using the linear transformation
$(e_1,\alpha,\beta)\longmapsto(e_1,\bar{\alpha},\beta)$ one may
write $\mathcal{C}=(\bar{\alpha}, \widetilde{\gamma}, \bar{\alpha},
x_4,...,x_m)$. Furthermore,  we can carry out our argument as in
Proposition~\ref{c2} to reduce $\mathcal{C}$ to $(\bar{\alpha},
\widetilde{\gamma}, \bar{\alpha}, \beta, \alpha, \beta,..., \alpha,
\beta)$. If $\widetilde{\gamma}$ is not $\gamma$, applying the
operation $\text{O}_{22}$ on $\alpha,\bar{\gamma},\bar{\alpha}$, one
may further reduce $\mathcal{C}'$ to $\mathcal{C}_8$.
\end{proof}

\begin{rem}\label{c4}
An easy observation shows that for a 3-independent trivial coloring
$\lambda$ with the property $(\star)$,  if $m=3$, then $\lambda$ is
just sector-equivalent to the following canonical form
$$\lambda_{C^3} \text{ with the coloring sequence }
\mathcal{C}^3=(\bar{\gamma},\alpha,\beta)$$ and if $m=4$, then
$\lambda$ is just sector-equivalent to one of  the following two
canonical forms
\begin{enumerate}
\item  $\lambda_{C_1^4}$ with  the coloring sequence $\mathcal{C}_1^4=(\alpha,\gamma,\bar{\alpha},\beta).$
\item  $\lambda_{C_2^4}$ with  the coloring sequence $\mathcal{C}_2^4=(\alpha,\bar{\gamma},\alpha,\beta).$
\end{enumerate}
\end{rem}

Combining Propositions~\ref{c1}-\ref{c3} and Remark~\ref{c4} gives
the following

\begin{cor}\label{nu}
The number  of homeomorphism classes of small covers over
$\Pol^3(m)$ with  trivial colorings is at most
$$N_t(m)=\begin{cases}
2 & \text{ if } m=3\\
4  & \text{ if  $m>3$ is odd}\\
6  &\text{ if  $m$ is even} \end{cases}$$
\end{cor}

By Proposition~\ref{orient}, a direct observation shows that

\begin{cor}\label{ori-1}
 $M(\lambda_{C_i}), i=1,5, 10$, are orientable, and $M(\lambda_{C_i}),
i=2,3,4,6,7,8,9$, are non-orientable.
\end{cor}

\subsection{Nontrivial colorings} Given a pair $(\Pol^3(m), \lambda)$ in
$\Lambda(\Pol^3(m))$, throughout the following suppose that $\lambda$ is
nontrivial, i.e., $\lambda(c)\not=\lambda(f)$. \vskip .2cm

\begin{defn} Let  $m_\lambda$ denote the number of side faces in
$$M_\lambda=\{s_i \big|
 \text{Span}\{\lambda(s_{i-1}), \lambda(s_i),\lambda(s_{i+1}) \}=(\Z_2)^3 \}.$$ Set $\lambda_0:=\lambda(c)-\lambda(f)$. Let
$n_\lambda$ denote the number of side faces in
$N_\lambda=\{s_i\big|\lambda(s_i)=\lambda_0\}$.
\end{defn}

\begin{lem}\label{lemnontri}
Let $\lambda$ be a nontrivial coloring on $\Pol^3(m)$ with $m>3$.
Then
\begin{enumerate}
\item $m_\lambda\leq n_\lambda\leq m/2$. In particular, if $m$ is odd, then $n_\lambda>0$.
\item $m_\lambda$ is  even.
\end{enumerate}
\end{lem}

\begin{proof} First, $n_\lambda\leq m/2$ is obvious since any two
faces in $N_\lambda$ are not  adjacent. To show that $m_\lambda\leq
n_\lambda$, take one  $s_i\in N_3$. Then
$\lambda(s_{i-1}),\lambda(s_i),\lambda(s_{i+1})$ are linearly
independent.  Furthermore, the linear independence of
$\{\lambda(s_{i-1}),\lambda(s_i),\lambda(c)\}$ and
$\{\lambda(s_i),\lambda(s_{i+1}),\lambda(c)\}$ implies that
$\lambda(c)$ must be either $\lambda(s_{i-1})+\lambda(s_{i+1})$ or
$\lambda(s_{i-1})+\lambda(s_i)+\lambda(s_{i+1})$. This is also true
for $\lambda(f)$. Now $\lambda(c)\neq\lambda(f)$ makes sure that
$\lambda(c)-\lambda(f)=\lambda(s_i)$, so that $s_i\in N_\lambda$.
Thus,  $M_\lambda\subseteq N_\lambda$, i.e., $m_\lambda\leq
n_\lambda$. Moreover, if $n_\lambda=0$ then $m_\lambda=0$, so that
the coloring sequence of $\lambda$ is 2-independent. However,
$n_\lambda=0$ means that one has only two choices of colors. This
forces $m$ to be even.

With no loss, assume that $m_\lambda>0$ (since $0$ is even). For
each $i$, let $V_i$ denote the subspace spanned by $\lambda(s_i)$
and $\lambda(s_{i+1})$.  Obviously, if $s_i\in M_\lambda$ then  $
V_i\not= V_{i-1}$, and if $s_i\not\in M_\lambda$ then $V_i=V_{i-1}$.
Thus, $s_i\in M_\lambda$ if and only if $ V_i\not= V_{i-1}$. Next we
claim that for any $i$, $\lambda_0\in V_i$. In fact, if $s_i\in
M_\lambda$, since $M_\lambda\subseteq N_\lambda$, we have that
$s_i\in N_\lambda$ so $\lambda_0=\lambda(s_i)\in V_i$. If
$s_i\not\in M_\lambda$ then $V_i=V_{i-1}$. Since  $V_i=V_{i-1}$
contains no $\lambda(c)$ and $\lambda(f)$, $\lambda_0$ must be in
$V_i=V_{i-1}$. This proves the claim. Furthermore, for any $i$,
$V_i$ must be either $\text{Span}\{\lambda_0, \alpha\}$ or
$\text{Span}\{\lambda_0, \lambda(c)+\alpha\}$, where $\alpha$ is a
nonzero element such that $\alpha$, $\lambda(c)$ and $\lambda(f)$
are linearly independent. Now we clearly see that there is a switch
of choosing either $\text{Span}\{\lambda_0, \alpha\}$ or
$\text{Span}\{\lambda_0, \lambda(c)+\alpha\}$ exactly when we pass
$s_i\in M_\lambda$. But the total number of switches must be even.
So $m_\lambda$ is even.
\end{proof}

\begin{rem}
An easy observation shows that if $m=3$, then there is a possibility
that $m_\lambda=3$ but still $n_\lambda=1<3/2$. This is exactly an
exception only for $m_\lambda$ in the case $m=3$.
\end{rem}

Throughout the following, assume that $m>3$.

\vskip .2cm

Applying sectors $S(2_1), S(3_1), S(4)$ and $S(5)$ to the nontrivial
coloring $\lambda$ gives the following four fundamental operations
on its coloring sequence:

\begin{enumerate}
\item[$\bar{\text{O}}_{21}$] Take $s_k, s_l\in N_\lambda$ with $k<l$, by using $S(2_1)$, we may do a linear
transformation $(\lambda(c), \lambda_0,
\lambda(s_{k+1})\mapsto(\lambda(c), \lambda_0,
\lambda(s_{k+1})+\lambda_0)$ to change the coloring sequence of
$s_k, s_{k+1}, ..., s_l$.
\item[$\text{O}_{31}$] Take $s_k, s_l$ $(k<l)$ with $\lambda(s_k)=\lambda(s_l)\neq \lambda_0$, and use $S(3_1)$ to reflect the
coloring sequence of $s_k, s_{k+1}, ..., s_l$.
\item[$\text{O}_{4}$] Take $s_k, s_l$ $(k<l)$ with $\lambda(s_k)-\lambda(s_l)=\lambda(c)$ or $\lambda(f)$,  we
use $S(4)$ to reflect the coloring sequence of $s_k, s_{k+1}, ...,
s_l$ and then to do a linear transformation
$(\lambda_0,\lambda(s_k),\lambda(s_l))\mapsto(\lambda_0,\lambda(s_l),\lambda(s_k))$
to change the reflected coloring sequence.
\item[$\text{O}_{5}$] Take $s_k, s_l$ $(k<l)$ with
$\lambda(s_k)-\lambda(s_l)=\lambda_0$, we use $S(4)$ to reflect the
coloring sequence of $s_k, s_{k+1}, ..., s_l$ and then to do a
linear transformation $(\lambda(c),\lambda(s_k),$ $\lambda(s_l))$
$\mapsto(\lambda(c),\lambda(s_l),\lambda(s_k))$ to change the
reflected coloring sequence.
\end{enumerate}

It is easy to check the following

\begin{lem}\label{op}
The operations $\bar{\text{\rm O}}_{21}$, $\text{\rm O}_{31}$,
$\text{\rm O}_{4}$ and $\text{\rm O}_{5}$ above will not change $
m_\lambda, n_\lambda$ of the nontrivial coloring $\lambda$.
\end{lem}

 Without loss of generality, throughout the following  one assumes that $\lambda(c)=e_1,\lambda(f)=e_1+e_2$, so
that $\lambda_0=e_2$, where  $\{e_1,e_2,e_3\}$ is a basis of
$(\Z_2)^3$.

\begin{prop}\label{pNONTRI}
For $m>3$, each nontrivial coloring $\lambda$ is sector-equivalent
to the following canonical form $\lambda_{C_*}$ with the coloring
sequence
\begin{equation}\label{c*}
\mathcal{C}_*=(e_2, x_1, e_2, ..., e_2, x_{m_\lambda}, e_2, y_1,
..., e_2, y_{n_\lambda-m_\lambda},  z_1, ...,
z_{m-2n_\lambda})\end{equation}
where $x_i=\begin{cases} e_1+e_3 & \text{if $i$ is odd}\\
e_3 & \text{if $i$ is even} \end{cases}$ and for all $1\leq i\leq
n_\lambda-m_\lambda$, $y_i=e_3$
and $z_i=\begin{cases} e_2+e_3 & \text{if $i$ is odd}\\
e_3 & \text{if $i$ is even.} \end{cases}$
\end{prop}

\begin{proof}
If $n_\lambda\leq 1$ then clearly the coloring $\lambda$ can be
reduced to a coloring with the coloring sequence
$$\begin{cases} (e_2, e_3, e_2+e_3, e_3, ..., e_2+e_3, e_3) & \text{ if $n_\lambda=1$ and $m$ is even}\\
(e_2, e_3, e_2+e_3, e_3, ..., e_2+e_3, e_3, e_2+e_3) & \text{ if $n_\lambda=1$ and $m$ is odd}\\
(e_2+e_3, e_3, ..., e_2+e_3, e_3). & \text{ if
$n_\lambda=0$}\end{cases}$$ If $n_\lambda\geq 2$, we may choose two
$s_k$ and $s_l$ in $N_\lambda$ with $k<l$. Consider the coloring
sub-sequence
\begin{equation}\label{ntc}
(\lambda(s_k)=)e_2, r_1, ..., r_{l-2},e_2(=\lambda(s_l)),r_{l-1}
\end{equation} of
$s_k,..., s_l, s_{l+1} $, it is easy to see that $r_1-r_{l-1}\in
\text{Span}\{e_1,e_2\}$. Then when  $r_1-r_{l-1}=0$ (resp. $e_1$
or $e_1+e_2$, $e_2$),  we may do the operation $\text{O}_{31}$
(resp. $\text{O}_4$ or $\text{O}_5$) on $r_1, ..., r_{l-2},e_2,
r_{l-1}$ from $s_{k+1}$ to $s_{l+1}$, and change (\ref{ntc}) into
$e_2, r_{l-1}, e_2, r'_{l-2},..., r'_2, r_1$. With this understood,
assume that $N_\lambda=\{s_1$, $s_3,$$ ...,$ $s_{2n_\lambda-1}\}$,
so we may write the coloring sequence of $\lambda$ as follows:
$$\mathcal{C}=(e_2, \alpha_1, ..., e_2, \alpha_{n_\lambda}, \beta_1,
..., \beta_{m-2n_\lambda})$$ with $\alpha_{n_\lambda}, \beta_i\in
\{e_2+e_3,e_3\}$. By doing the operation $\bar{\text{O}}_{21}$ on
$\alpha_{n_\lambda}, \beta_1, ..., \beta_{m-2n_\lambda}$, we may
reduce $\mathcal{C}$ to $\mathcal{C}'=(e_2, \alpha_1, ...,
\alpha_{n_\lambda-1}, e_2, e_3, z_1, ..., z_{m-2n_\lambda})$ such
that $z_i$ is $e_2+e_3$ if $i$ is odd, and $e_3$ if $i$ is even.
Then we may further use the operation $\text{O}_4$ to reduce
$\mathcal{C}'$ to $\mathcal{C}''$ with $M_\lambda=\{s_1, ...,
s_{2m_\lambda-1}\}$ and without changing the part
$\mathcal{C}'-\{e_2, \alpha_1, ..., \alpha_{n_\lambda-1}, e_2\}$.
Finally, by using the operation $\bar{\text{O}}_{21}$, we may reduce
$\mathcal{C}''$ to $\mathcal{C}_*$ as desired.
\end{proof}

Together with Theorem~\ref{t:GENERAL},
Lemmas~\ref{lemnontri}-\ref{op} and Proposition~\ref{pNONTRI}, it
easily follows that
\begin{cor}\label{home}
Let $\lambda_1, \lambda_2$ be two nontrivial colorings on
$\Pol^3(m)$ with $m>3$. If $(m_{\lambda_1},
n_{\lambda_1})=(m_{\lambda_2}, n_{\lambda_2})$, then $M(\lambda_1)$
and $M(\lambda_2)$ are homeomorphic.
\end{cor}

\begin{cor}\label{num}
For $m>3$, let $(k,l)$ be a pair such that $(1)$ $l\leq k\leq m/2$
and if $2\nmid m$ then $k>0$; and $(2)$ $l$ is even. Then there is a
nontrivial coloring $\lambda$ on $\Pol^3(m)$ with
$(n_\lambda,m_\lambda)=(k,l)$.
\end{cor}

As a consequence of Proposition~\ref{orient} and
Proposition~\ref{pNONTRI}, one also has
\begin{cor}
Let $\lambda$ be a nontrivial coloring on $\Pol^3(m)$ with $m>3$.
Then $M(\lambda)$ is orientable if $n_\lambda=0$, and non-orientable
if $n_\lambda>0$.
\end{cor}

\section{Mod 2 cohomology rings and two invariants} \label{cohomology}

Given a pair $(\Pol^3(m), \lambda)$ in $\Lambda(\Pol^3(m))$, one
knows that the mod 2 cohomology ring of $M(\lambda)$ is
$$H^*(M(\lambda);\Z_2)=\Z_2[c,f, s_1,...,s_m]/I+J_\lambda$$
where $I$ is the ideal generated by $cf$ and $s_is_j$ with $s_i\cap
s_j=\emptyset$, and $J_\lambda$ is the ideal generated by three
linear relations (determined by the $3\times (m+2)$ matrix
$(\lambda(c), \lambda(f), \lambda(s_1), ..., \lambda(s_m))$).

\subsection{Two invariants $\Delta(\lambda)$ and $\mathcal{B}(\lambda)$} Now let us introduce two invariants
in
 $H^*(M(\lambda)$; $\Z_2)$. Set
$$\mathcal{H}_\lambda^1=\{x\in H^1(M(\lambda); \Z_2)\big|x^2=0\}$$ and
$$\mathcal{H}_\lambda^2=\{x^2\big|x\in H^1(M(\lambda); \Z_2)\}.$$ Obviously,  they are all vector
spaces over $\Z_2$, and $\dim\mathcal{H}_\lambda^1$ is  an invariant
of the cohomology ring $H^*(M(\lambda);\Z_2)$, denoted by
$\Delta(\lambda)$.  So
$$\dim\mathcal{H}_\lambda^2=\dim
H^1(M(\lambda);
\Z_2)-\dim\mathcal{H}_\lambda^1=m-1-\Delta(\lambda).$$
Note that $\dim H^1(M(\lambda); \Z_2)=m-1$ by Example~\ref{betti}.

Consider the bilinear map $$\omega:H^1(M(\lambda);\Z_2)\times
\mathcal{H}_\lambda^1\longrightarrow H^2(M(\lambda);\Z_2)$$
by $(x, y)\longmapsto xy$. Let $\mathcal{K}_\lambda=\text{Im}\omega$. Then $\mathcal{K}_\lambda\cong H^1(M(\lambda);
\Z_2)\otimes \mathcal{H}_\lambda^1.$ Furthermore, we can induce  a bilinear map
$$\tilde{\omega}: H^1(M(\lambda);\Z_2)\times
\mathcal{H}_\lambda^1\longrightarrow
\mathcal{K}_\lambda/(\mathcal{K}_\lambda\cap\mathcal{H}_\lambda^2)$$
by $(x, y)\longmapsto [xy]$, which is surjective. Let
$\text{Hom}(\mathcal{K}_\lambda/(\mathcal{K}_\lambda\cap\mathcal{H}_\lambda^2),
\Z_2)$ be  the dual space of
$\mathcal{K}_\lambda/(\mathcal{K}_\lambda\cap\mathcal{H}_\lambda^2)$ as a vector space over $\Z_2$.
Take a
$\theta\in\text{Hom}(\mathcal{K}_\lambda/(\mathcal{K}_\lambda\cap\mathcal{H}_\lambda^2),
\Z_2)$, one can obtain a bilinear function $$\theta\circ
\tilde{\omega}: H^1(M(\lambda);\Z_2)\times
\mathcal{H}_\lambda^1\longrightarrow\Z_2.$$
Now given a basis $\{\alpha_1, ..., \alpha_{m-1}\}$ of $H^1(M(\lambda);\Z_2)$ and a basis $\{\beta_1, ..., \beta_{\Delta(\lambda)}\}$ of
$\mathcal{H}_\lambda^1$. Then $\theta\circ
\tilde{\omega}$  corresponds an
$(m-1)\times \Delta(\lambda)$-matrix $M_\theta=(m_{ij})$ where $m_{ij}=\theta([\alpha_i\beta_j])$.  It is easy to check that  the rank of $M_\theta$ does not depend upon
 the choices of the bases of $H^1(M(\lambda);\Z_2)$ and
$\mathcal{H}_\lambda^1$, but it depends upon the choice of $\theta$.  Let $b_r(\lambda)$
 denote the number of those
$\theta\in\text{Hom}(\mathcal{K}_\lambda/(\mathcal{K}_\lambda\cap\mathcal{H}_\lambda^2),
\Z_2)$ such that $\text{rank } M_\theta=r$ where $1\leq
r\leq \Delta(\lambda)$. Then we obtain an integer vector
$$\mathcal{B}(\lambda)=(b_1(\lambda), ...,
b_{\Delta(\lambda)}(\lambda)).$$
It is not difficult to see that $\mathcal{B}(\lambda)$ is
 an invariant of the cohomology ring $H^*(M(\lambda);\Z_2)$.

It should be pointed out that we shall only calculate $b_1(\lambda)$
and $b_2(\lambda)$ in $\mathcal{B}(\lambda)$ because  this
will  be sufficient enough to reach our purpose. By
$\bar{\mathcal{B}}(\lambda)$ we denote $(b_1(\lambda),
b_2(\lambda))$. We also use the convention that $b_2(\lambda)=0$ if
$\Delta(\lambda)=1$.

\subsection{Calculation of $\Delta(\lambda)$}
First we shall deal with the case in which $\lambda$ is nontrivial.

\begin{lem}\label{alg-1}
If $\lambda$ is nontrivial, then
$$\Delta(\lambda)=\begin{cases}
n_\lambda & \text{if $n_\lambda>0$ and  $m_\lambda=0$}\\
n_\lambda-1 & \text{if $m_\lambda>0$}\\
1 & \text{if $n_\lambda=0$ (so $m$ is even).}
\end{cases}$$
\end{lem}

\begin{proof}
By Proposition~\ref{pNONTRI}, each $\lambda$ is  sector-equivalent
to the canonical form $\lambda_{C_*}$ with the coloring sequence
$\mathcal{C}_*$ and without changing $m_\lambda$ and $n_\lambda$, so
it suffices to consider the  $\lambda_{C_*}$.

 If $n_\lambda>0$ and $m_\lambda=0$, we can obtain from
(\ref{c*}) that $\lambda_{C_*}$ determines the following three
linear relations in $H^1(M(\lambda_{C_*});\Z_2)$
\begin{equation}\label{r1}
c+f=0
\end{equation}
\begin{equation}\label{r2}
f+\sum_{i \text{ is odd}}s_i=0
\end{equation}
\begin{equation}\label{r3}
s_2+\cdots+s_{2n_\lambda}+\sum_{2n_\lambda <i\leq m}s_i=0.
\end{equation}
So we may choose $B_1=\{f, s_2, s_3, ..., s_{m-1}\}$  as a basis of
$H^1(M(\lambda_{C_*});\Z_2)$. Since $cf=0$ and $s_is_j=0$ with
$s_i\cap s_j=\emptyset$ in $H^*(M(\lambda_{C_*});\Z_2)$, one can
easily obtain from (\ref{r1}) and (\ref{r3}) that $B_2=\{f, s_2,
s_4, ..., s_{2n_\lambda-2}\}\subset \mathcal{H}_{\lambda_{C_*}}^1$,
and $B_2\subset B_1$. Thus,
$\dim\mathcal{H}_\lambda^1=\dim\mathcal{H}_{\lambda_{C_*}}^1\geq
n_\lambda$. On the other hand, an easy argument shows that
$B_3=\{s_3^2, ..., s_{2n_\lambda-1}^2, s_{2n_\lambda}^2, ...,
s_{m-1}^2, fs_2, s_3s_4, s_5s_6, ..., s_{2n_\lambda-1}s_{2n_\lambda}
\}$ forms a basis of $H^2(M(\lambda_{C_*});\Z_2)$. Now observe that
the square of each element of $B_1\setminus B_2$ is in $B_3$, so
$\dim\mathcal{H}_\lambda^2=\dim\mathcal{H}_{\lambda_{C_*}}^2\geq
m-1-n_\lambda$. Furthermore, $\dim\mathcal{H}_\lambda^1\leq
n_\lambda$. Therefore, $\Delta(\lambda)=n_\lambda$.

\vskip .2cm If $m_\lambda>0$, then $\lambda_{C_*}$ determines the
following three linear relations
$$\begin{cases}
 c+f+s_2+\cdots+s_{2m_\lambda-2}=0 \\
f+\sum_{i \text{ is odd}}s_i=0\\
s_2+\cdots+s_{2n_\lambda}+\sum_{2n_\lambda <i\leq m}s_i=0.
\end{cases}$$
In this case, we choose $B_4=\{s_1, s_2, ..., s_{m-1}\}$ as a basis
of $H^1(M(\lambda_{C_*});\Z_2)$. Then one sees that $B_5=\{s_2, s_4,
..., s_{2n_\lambda-2}\}\subset \mathcal{H}_{\lambda_{C_*}}^1$.
Furthermore, we choose $$B_6=\{s_1^2, s_3^2, ...,
s_{2n_\lambda-1}^2, s_{2n_\lambda}^2, s_{2n_\lambda+1}^2, ...,
s_{m-1}^2, s_2s_3, s_4s_5, ..., s_{2n_\lambda-2}s_{2n_\lambda-1}\}$$
as a basis of $H^2(M(\lambda_{C_*});\Z_2)$. A similar argument as
above shows that $\Delta(\lambda)=n_\lambda-1$.

\vskip .2cm If $n_\lambda=0$, in a similar way as above, it is easy
to see that we may choose $B_7=\{f, s_3,  ..., s_m\}$ as a basis of
$H^1(M(\lambda_{C_*});\Z_2)$ and  $B_8=\{f\}$ forms a basis of
$\mathcal{H}_{\lambda_{C_*}}^1$. Thus, $\Delta(\lambda)=1$.
\end{proof}

Next we consider the case in which $\lambda$ is trivial.

\begin{lem}\label{alg-2}
Let $\lambda$ be trivial. Then
$$\Delta(\lambda)=\begin{cases}
m-1 & \text{if $\lambda\approx\lambda_{C_1}$}\\
m-2 & \text{if $\lambda\approx\lambda_{C_i}, i=2,3,4$}\\
m-3 & \text{if $\lambda\approx\lambda_{C_i}$ with $m>4$,
$i=5,6,7,8,9, 10$.}
\end{cases}$$
In particular, if $m=3$ then $\Delta(\lambda_{C^3})=0$, and if $m=4$
then $\Delta(\lambda_{C_1^4})=\Delta(\lambda_{C_2^4})=1$.
\end{lem}
\begin{proof}
The argument is similar to that of Lemma~\ref{alg-1}, and is not
quite difficult. Here we only list the three linear relations and
the bases of $H^i(M(\lambda);\Z_2) (i=1,2)$ and
$\mathcal{H}_\lambda^1$, but for the detailed proof, we would like
to leave it to readers as an exercise.
\begin{center}
\begin{tabular}{|c|c|c|}
\hline $\lambda$ & $m$ & Three linear relations by determined by $J_\lambda$ \\
\hline
   $\lambda_{C_1}$  & even  & $c+f=0, \sum_{i\text{ is odd}}s_i=0, \sum_{i\text{ is even}}s_i=0$     \\
\hline
    $\lambda_{C_2}$ & even  &    $c+f=0, s_2+\sum_{i\text{ is odd}}s_i=0, \sum_{i\text{ is even}}s_i=0$  \\
\hline
  $\lambda_{C_3}$     & odd & $c+f=0, s_2+\sum_{i\text{ is odd}}s_i=0, \sum_{i\text{ is even}}s_i=0$      \\
\hline
    $\lambda_{C_4}$     & even  & $c+f+s_1=0, \sum_{i\text{ is odd}}s_i=0, \sum_{i\text{ is even}}s_i=0$  \\
\hline
  $\lambda_{C_5}$  & odd  & $c+f+s_1=0, \sum_{i\text{ is odd}}s_i=0, s_1+\sum_{i\text{ is even}}s_i=0$  \\
\hline
 $\lambda_{C_6}$  & odd  & $c+f+s_1+s_2=0, \sum_{i\text{ is odd}}s_i=0, s_1+\sum_{i\text{ is even}}s_i=0$  \\
\hline
 $\lambda_{C_7}$  & odd  & $c+f+s_2=0, \sum_{i\text{ is odd}}s_i=0, s_1+\sum_{i\text{ is even}}s_i=0$  \\
\hline
 $\lambda_{C_8}$  & even  & $c+f+s_1+s_3=0, s_2+\sum_{i\text{ is odd}}s_i=0, \sum_{i\text{ is even}}s_i=0$  \\
\hline
 $\lambda_{C_9}$  & even   & $c+f+s_3=0, s_2+\sum_{i\text{ is odd}}s_i=0, \sum_{i\text{ is even}}s_i=0$  \\
\hline
 $\lambda_{C_{10}}$  & even   & $c+f+s_2=0, s_2+\sum_{i\text{ is odd}}s_i=0, \sum_{i\text{ is even}}s_i=0$  \\
\hline
\end{tabular}
\end{center}

\begin{center}
\begin{tabular}{|c|c|c|c|}
\hline $\lambda$ & Basis of $H^1(M(\lambda);\Z_2)$ & Basis of $\mathcal{H}_\lambda^1$ &  Basis of $H^2(M(\lambda);\Z_2)$\\
\hline
   $\lambda_{C_1}$  & $\{f, s_3, ..., s_m\}$  & $\{f, s_3, ..., s_m\}$ &$\{s_3s_4, fs_3, ...,fs_m\}$    \\
\hline
    $\lambda_{C_2}$ & $\{f, s_2, ..., s_{m-1}\}$  & $\{f, s_2, s_4, ..., s_{m-1}\}$ &  $\{s_1^2, fs_2, ..., fs_{m-1}\}$   \\
\hline
  $\lambda_{C_3}$     & $\{f, s_2, ..., s_{m-1}\}$  &   $\{f, s_4, ..., s_{m}\}$  &  $\{s_1^2, fs_1, ..., fs_{m-2}\}$ \\
\hline
    $\lambda_{C_4}$     & $\{f, s_2, ..., s_{m-1}\}$  &  $\{ s_3, ..., s_{m}\}$ &  $\{s_1s_2, fs_3, ..., fs_{m}\}$\\
\hline
  $\lambda_{C_5}$  & $\{f, s_2, ..., s_{m-1}\}$  &  $\{ s_3, ..., s_{m-1}\}$ &  $\{s_1s_2, fs_2, ..., fs_{m-1}\}$\\
\hline
 $\lambda_{C_6}$  & $\{f, s_2, ..., s_{m-1}\}$  &  $\{ s_3, ..., s_{m-1}\}$ &  $\{f^2, s_2^2, fs_2, ..., fs_{m-2}\}$ \\
\hline
 $\lambda_{C_7}$  & $\{f, s_2, ..., s_{m-1}\}$ &  $\{ s_3, ..., s_{m-1}\}$ &  $\{f^2, s_2^2, fs_3, ..., fs_{m-1}\}$ \\
\hline
 $\lambda_{C_8}$  & $\{f, s_2, ..., s_{m-1}\}$ & $\{ s_2, s_4, ..., s_{m-1}\}$ &  $\{f^2, s_3^2, fs_1, fs_4, ..., fs_{m-1}\}$ \\
\hline
 $\lambda_{C_9}$  &$\{f, s_2, ..., s_{m-1}\}$ &  $\{ s_2, s_4, ..., s_{m-1}\}$ &  $\{f^2, s_3^2, fs_1, fs_4,..., fs_{m-1}\}$ \\
\hline
 $\lambda_{C_{10}}$  & $\{f, s_2, ..., s_{m-1}\}$ &  $\{ s_2, s_4, ..., s_{m-1}\}$ &  $\{f^2, s_3^2, fs_3, ..., fs_{m-1}\}$\\
\hline
\end{tabular}
\end{center}
\end{proof}

\begin{rem}
Although it is not mentioned in this paper, the authors have
calculated the first Betti number under $\Z$-coefficients of all
small covers over prisms and discovered that the number is always
equal to $\Delta(\lambda)$ in the $\Z_2$-cohomology ring. One can
check that this is also true for all closed surfaces (i.e.,
2-dimensional small covers). It should be reasonable to conjecture
that this is true for all small covers.
\end{rem}

\begin{prop}\label{m>6}
Let $\lambda_1$ and $\lambda_2$ be two colorings in
$\Lambda(\Pol^3(m))$ such that $\lambda_1$ is trivial but
$\lambda_2$ is nontrivial. If $m>6$, then both $M(\lambda_1)$ and
$M(\lambda_1)$ cannot be homeomorphic.
\end{prop}
\begin{proof}
Suppose that $M(\lambda_1)$ and $M(\lambda_1)$ are homeomorphic.
Then their cohomologies are isomorphic, so
$\Delta(\lambda_1)=\Delta(\lambda_2)$. However, by
Lemmas~\ref{alg-1} and \ref{alg-2}, one has that
$\Delta(\lambda_1)\geq m-3$ and $\Delta(\lambda_2)\leq m/2$.
Furthermore, if $m>6$, then $\Delta(\lambda_1)\geq m-3>m/2\geq
\Delta(\lambda_2)$, so $\Delta(\lambda_1)\not=\Delta(\lambda_2)$, a
contradiction.
\end{proof}

\begin{rem}\label{3,4}
We see from the proof of Proposition~\ref{m>6} that
$\Delta(\lambda_1)$ and $\Delta(\lambda_2)$ can coincide only if
$m\leq 6$. For $m=5,6$, all possible cases that
$\Delta(\lambda_1)=\Delta(\lambda_2)$ happens are stated as follows:
 when $(n_{\lambda_{C_*}}, m_{\lambda_{C_*}})=(m-3, 0)$, one has that $\Delta(\lambda_{C_*})=\Delta(\lambda_{C_i})$, $i=5,
6, 7, 8, 9, 10$. For $m=3,4$, we know from \cite{ly} and \cite{cmo}
that up to homeomorphism,  there are only two small covers over
$\Pol^3(3)$: ${\Bbb R}P^3\sharp{\Bbb R}P^3$ and $S^1\times {\Bbb R}P^2$, and  there
are only four small covers over $\Pol^3(4)$: $(S^1)^3$, $S^1\times
K$,  a twist $(S^1)^2$-bundle over $S^1$ and a twist $K$-bundle over
$S^1$, where $K$ is a Klein bottle. In particular, the cohomological
rigidity holds in these cases.
\end{rem}

\subsection{Calculation of $\bar{\mathcal{B}}(\lambda)$} Let
$\lambda\in\Lambda(\Pol^3(m))$ with $m>4$. Choose an ordered basis
$B'$ of $H^1(M(\lambda);\Z_2)$ and an ordered basis $B''$ of
$\mathcal{H}_\lambda^1$, let  $A_0$ denote an
$(m-1)\times\Delta(\lambda)$ matrix $(a_{ij})$, where
$a_{ij}=[u_iv_j]\in
\mathcal{K}_\lambda/(\mathcal{K}_\lambda\cap\mathcal{H}_\lambda^2)$,
$u_i$ is the $i$-th element in $B'$ and $v_j$ the $j$-th element in
$B''$, so each element in $B'$ corresponds to a row and each element
in $B''$ a column. It follows that for any $\theta\in
\text{Hom}(\mathcal{K}_\lambda/(\mathcal{K}_\lambda\cap\mathcal{H}_\lambda^2),\Z_2)$,
$\theta(A_0)=(\theta(a_{ij}))$ is a representation matrix of
$\theta\circ\tilde{\omega}$.

\vskip .2cm
 First let us look at the case in which $\lambda$ is
nontrivial.
\begin{lem}\label{alg-3}
Let $\lambda$ be nontrivial. Then
$$\bar{\mathcal{B}}(\lambda)=\begin{cases}
(0,0) & \text{if $(n_\lambda, m_\lambda)=(0,0)$}\\
(1,0) & \text{if $(n_\lambda, m_\lambda)=(1,0)$ or $(2,2)$}\\
(1,3) & \text{if $(n_\lambda, m_\lambda)=(2,0)$}\\
(0,n_\lambda) & \text{if $n_\lambda>2$ and $m_\lambda=0$}\\
(n_\lambda-m_\lambda,{{m_\lambda-1}\choose 1}+{{m_\lambda-1}\choose
2}+{{n_\lambda-m_\lambda}\choose 2}) & \text{if $n_\lambda>2$ and
$m_\lambda>0$}
\end{cases}$$
\end{lem}
\begin{proof}
By Proposition~\ref{pNONTRI}, one may assume that
$\lambda=\lambda_{C_*}$. Then our argument proceeds as follows.

\vskip .2cm

(1) If $n_\lambda>0$ and $m_\lambda=0$, then from the proof of Lemma~\ref{alg-1} we
may take $B'=B_1$ and $B''=B_2$. Thus one has that $A_0$ is equal to
{\small \begin{equation*}
\begin{bmatrix}
0 &   [fs_2]    &  [fs_4] & [fs_6] & \cdots & [fs_{2n_\lambda-6}] & [fs_{2n_\lambda-4}] & [fs_{2n_\lambda-2}] \\
[s_2f]  & 0 &  0  & 0&  \hdotsfor{3}   &0 \\
  [s_3f]   & [s_3s_2] & [s_3s_4] & 0& \hdotsfor{3} & 0\\
  [s_4f] & 0 & 0& 0& \hdotsfor{3}   &0 \\
       \hdotsfor{8}          \\
  [s_{2n_\lambda-3}f] & 0 & \hdotsfor{3}& 0  & [s_{2n_\lambda-3}s_{2n_\lambda-4}]  & [s_{2n_\lambda-3}s_{2n_\lambda-2}]\\
   [s_{2n_\lambda-2}f]  & 0 &  0  & 0&  \hdotsfor{3}   &0 \\
[s_{2n_\lambda-1}f]  & 0 &  0  & 0&  \hdotsfor{2}   &0 & [s_{2n_\lambda-1}s_{2n_\lambda-2}]\\
[s_{2n_\lambda}f]  & 0 &  0  & 0&  \hdotsfor{3}   &0 \\
   \hdotsfor{8} \\
   [s_{m-1}f]  & 0 &  0  & 0&  \hdotsfor{3}   &0
\end{bmatrix}
\end{equation*}}

\noindent By direct calculations one knows from (\ref{r2}) and
(\ref{r3}) that
$s_{2n_\lambda}s_{2n_\lambda+1}+s_{2n_\lambda+1}^2+\cdots+s_m^2=0$
so $[s_{2n_\lambda}s_{2n_\lambda+1}]=0$ and
$$\begin{cases} s_{2i}s_{2i+1}=s_{2i+1}s_{2i+2}
& \text{when $1\leq i\leq n_\lambda-1$}\\
fs_i=s_i^2 \text{ so $[fs_i]=0$} & \text{when either $i$ is odd or $i>2n_\lambda$ is even}\\
fs_i=s_{i-1}s_i+s_is_{i+1} & \text{when $2\leq i\leq 2n_\lambda$ is
even.}
\end{cases}$$
Set $x_1=[fs_2]$ and $x_i=[s_{2i-1}s_{2i}]$ for $2\leq i\leq
n_\lambda$. Then $[fs_{2i}]=x_i+x_{i+1}$ for $2\leq i\leq
n_\lambda-1$ and $[fs_{2n_\lambda}]=x_{n_\lambda}$. Thus, we see
that $\{x_1,  ..., x_{n_\lambda}\}$ forms a basis of
$\mathcal{K}_\lambda/(\mathcal{K}_\lambda\cap\mathcal{H}_\lambda^2)$,
and the corresponding rows of $f, s_2, ..., s_{2n_\lambda}$ in $A_0$
are nonzero. Now we may reduce $A_0$ to $A$ by deleting those zero
rows of $A_0$, so that for each $\theta\in
\text{Hom}(\mathcal{K}_\lambda/(\mathcal{K}_\lambda\cap\mathcal{H}_\lambda^2),\Z_2)$,
$\text{rank}(\theta(A_0))=\text{rank}(\theta(A))$ still holds. Write
$A$ as follows: {\tiny \begin{equation*}
\begin{bmatrix}
0 &   x_1    &  x_2+x_3 & x_3+x_4 & \cdots & x_{n_\lambda-3}+x_{n_\lambda-2} & x_{n_\lambda-2}+x_{n_\lambda-1} & x_{n_\lambda-1}+x_{n_\lambda} \\
x_1  & 0 &  0  & 0&  \hdotsfor{3}   &0 \\
  0   & x_2 & x_2 & 0& \hdotsfor{3} & 0\\
   x_2+x_3 & 0 & 0& 0& \hdotsfor{3}   &0 \\
       \hdotsfor{8}          \\
  0 & 0 & \hdotsfor{3}& 0  & x_{n_\lambda-1}  & x_{n_\lambda-1}\\
   x_{n_\lambda-1}+x_{n_\lambda}  & 0 &  0  & 0&  \hdotsfor{3}   &0 \\
0  & 0 &  0  & 0&  \hdotsfor{2}   &0 & x_{n_\lambda}\\
x_{n_\lambda}  & 0 &  0  & 0&  \hdotsfor{3}   &0 \\
  \end{bmatrix}
\end{equation*}}

\noindent Let $\{\theta_i\big| i=1,..., n_\lambda\}$ be the dual
basis of $\{x_1, ..., x_{n_\lambda}\}$ in
$\text{Hom}(\mathcal{K}_\lambda/(\mathcal{K}_\lambda\cap\mathcal{H}_\lambda^2),\Z_2)$.
Take any $\theta\in
\text{Hom}(\mathcal{K}_\lambda/(\mathcal{K}_\lambda\cap\mathcal{H}_\lambda^2),\Z_2)$,
one may write $\theta=\sum_{i\in S}\theta_i$ where
$S\subset\{1,...,n_\lambda\}$. Obviously,  if $n_\lambda=1$ then
$b_1(\lambda)=1$ and $b_2(\lambda)=0$. If $n_\lambda\geq 2$, it is
easy to see that $b_1(\lambda)=0$ since $\text{rank}\theta(A)$
cannot be 1 whenever $S$ is empty or non-empty. If $n_\lambda=2$,
then $\text{rank}\theta(A)=2$ only when $S=\{1\}, \{2\}, \{1, 2\}$,
so $b_2(\lambda)=3$. If $n_\lambda>2$, by direct calculations, one
has that only when $S=\{i\} (i\not=2)$ or $\{1,2\}$,
$\text{rank}\theta(A)=2$, so $b_2(\lambda)=n_\lambda$.
 \vskip .2cm

(2) If $m_\lambda>0$, then from the proof of Lemma~\ref{alg-1} we may take $B'=B_4$
and $B''=B_5$. Moreover, we see that
 in $A_0$, only corresponding rows of
$s_1,s_3,...,s_{2n_\lambda-1}$ are nonzero, so we may delete the
other rows from $A_0$ to obtain $A$ so that for each $\theta\in
\text{Hom}(\mathcal{K}_\lambda/(\mathcal{K}_\lambda\cap\mathcal{H}_\lambda^2),\Z_2)$,
$\text{rank}(\theta(A_0))=\text{rank}(\theta(A))$. Now we can write
down $A$ after simple calculations:
\begin{equation*}
A=
\begin{bmatrix}
[s_1s_2] &   0    &  \hdotsfor{2}    & 0 \\
[s_2s_3] & [s_2s_3] &  0  &  \hdotsfor{1}   &0 \\
   0   & [s_4s_5] & [s_4s_5] & 0& \hdotsfor{1} \\
       \hdotsfor{5}          \\
       \hdotsfor{2}& 0  & [s_{2n_\lambda-4}s_{2n_\lambda-3}]  & [s_{2n_\lambda-4}s_{2n_\lambda-3}]\\
       \hdotsfor{3}& 0  & [s_{2n_\lambda-2}s_{2n_\lambda-1}]
\end{bmatrix}
\end{equation*}
Set $x_i=[s_{2i}s_{2i+1}]\in
\mathcal{K}_\lambda/(\mathcal{K}_\lambda\cap\mathcal{H}_\lambda^2)$,
for $i=1, ..., n_\lambda-1$. A direct calculation shows that
$[s_1s_2]=x_1+x_2+\cdots+x_{m_\lambda-1}$. So we see that
$\{x_i\big|i=1, ..., n_\lambda-1\}$ forms a basis of
$\mathcal{K}_\lambda/(\mathcal{K}_\lambda\cap\mathcal{H}_\lambda^2)$.
Let $\{\theta_i\big|i=1, ..., n_\lambda-1\}$ be its dual basis in
$\text{Hom}(\mathcal{K}_\lambda/(\mathcal{K}_\lambda\cap\mathcal{H}_\lambda^2),\Z_2)$.
Then one may write $\theta=\sum_{i\in S}\theta_i$ where
$S\subset\{1,...,n_\lambda-1\}$. Now $i\in S$ implies that the
$(i+1)$-th row of $\theta(A)$ is nonzero.  In order that
$\text{rank}(\theta(A))=1$, one must have $\sharp(S)=1$ since $S$
cannot be empty. If $n_\lambda=2$ then $m_\lambda=2$,
$\mathcal{H}^1_\lambda$ is 1-dimensional and
$\mathcal{K}_\lambda/(\mathcal{K}_\lambda\cap\mathcal{H}_\lambda^2)$
has only a nonzero element, so $b_1(\lambda)=1$ and
$b_2(\lambda)=0$. If $n_\lambda>2$, then $\text{rank}(\theta_i(A))=1$ if
and only if $\theta_i([s_1s_2])=0$. This is equivalent to that
$\theta_i(x_1+x_2+...+x_{m_\lambda-1})=0\Leftrightarrow
i>m_\lambda-1$. Therefore, $b_1(\lambda)=n_\lambda-m_\lambda$. In
this case, an easy argument shows that
$b_2(\lambda)={{m_\lambda-1}\choose 1}+{{m_\lambda-1}\choose
2}+{{n_\lambda-m_\lambda}\choose 2}$.

\vskip .2cm

(3) If $n_\lambda=0$, then from the proof of Lemma~\ref{alg-1} we may take $B'=B_4$
and $B''=B_5$, so
$$A_0=(0, [s_3f], ..., [s_mf]).$$
However, a direct calculation shows that for each $i$, $s_if=s_i^2$,
so $[s_if]=0$ in
$\mathcal{K}_\lambda/(\mathcal{K}_\lambda\cap\mathcal{H}_\lambda^2)$.
Thus,
$\dim\mathcal{K}_\lambda/(\mathcal{K}_\lambda\cap\mathcal{H}_\lambda^2)=0$,
and so  $\bar{\mathcal{B}}(\lambda)=(0,0)$.
\end{proof}

\begin{thm}\label{c-r1}
Let $\lambda_1, \lambda_2$ be two nontrivial colorings on
$\Pol^3(m)$ with $m>4$. Then $M(\lambda_1)$ and $M(\lambda_2)$ are
homeomorphic if and only if  their cohomologies
$H^*(M(\lambda_1);\Z_2)$ and $H^*(M(\lambda_2);\Z_2)$ are isomorphic
as rings.
\end{thm}
\begin{proof}
It suffices to show that if $H^*(M(\lambda_1);\Z_2)$ and
$H^*(M(\lambda_2);\Z_2)$ are isomorphic, then $M(\lambda_1)$ and
$M(\lambda_2)$ are homeomorphic. Now suppose that
$H^*(M(\lambda_1);\Z_2)\cong H^*(M(\lambda_2);\Z_2)$. Then  one has
that $\bar{\mathcal{B}}(\lambda_1)=\bar{\mathcal{B}}(\lambda_2)$. We
claim that $(m_{\lambda_1}, n_{\lambda_1})=(m_{\lambda_2},
n_{\lambda_2})$. If not, then by Lemma~\ref{alg-3}, the possible
case in which this happens is
$\bar{\mathcal{B}}(\lambda_1)=\bar{\mathcal{B}}(\lambda_2)=(1,0)$.
Without loss of generality, assume that $(m_{\lambda_1},
n_{\lambda_1})=(1,0)$ and $(m_{\lambda_2}, n_{\lambda_2})=(2,2)$.
Then by Lemma~\ref{alg-1}, one has
$\Delta(\lambda_1)=\Delta(\lambda_2)=1$, so
$\mathcal{H}^1_{\lambda_1}$ and $\mathcal{H}^1_{\lambda_2}$ contains
only a nonzero element. Let $z_0^{(i)}$ be the unique nonzero
element of $\mathcal{H}^1_{\lambda_i}, i=1,2$. For each $i$, define
a linear map $\Phi_i: H^1(M(\lambda_i);\Z_2)\longrightarrow
H^2(M(\lambda_i);\Z_2)$ by $x\longmapsto z_0^{(i)}x$.

\vskip .2cm When $i=1$, by Lemma~\ref{alg-1} one may choose
$B_1=\{f, s_2, s_3, ..., s_{m-1}\}$ as a basis of
$H^1(M(\lambda_1);\Z_2)$ and $B_2=\{f\}$ as a basis of
$\mathcal{H}^1_{\lambda_1}$, so $z_0^{(1)}=f$. By direct
calculations, one has that for $3\leq j\leq m-1$, $fs_j=s_j^2$.
Since $fs_2, s_3^2, ..., s_{m-1}^2$ are linearly independent, one
knows that $\Phi_1$ has rank $m-2$.

\vskip .2cm

When $i=2$, by Lemma~\ref{alg-1} one may choose $B_4=\{s_1, s_2,
..., s_{m-1}\}$ as a basis of $H^1(M(\lambda_2);\Z_2)$ and
$B_5=\{s_2\}$ as a basis of $\mathcal{H}^1_{\lambda_2}$, so
$z_0^{(1)}=s_2$. Since $s_2^2=s_2s_j=0, j\geq 4$, one sees that
$\Phi_2$ has rank at most $2$.

\vskip .2cm

Now since $m>4$, one has that
$\text{rank}\Phi_1=m-2>2\geq\text{rank}\Phi_2$, but this is
impossible. Thus, one must have $(m_{\lambda_1},
n_{\lambda_1})=(m_{\lambda_2}, n_{\lambda_2})$. Moreover, the
theorem follows from Corollary~\ref{home}.
\end{proof}

\begin{cor}\label{complete invariant}
Let $\lambda_1, \lambda_2$ be two nontrivial colorings on
$\Pol^3(m)$ with $m>4$. Then $M(\lambda_1)$ and $M(\lambda_2)$ are
homeomorphic if and only if $(m_{\lambda_1},
n_{\lambda_1})=(m_{\lambda_2}, n_{\lambda_2})$.
\end{cor}

Furthermore, by Corollary~\ref{num} one has

\begin{cor}\label{number-1}
The number  of homeomorphism classes of small covers over
$\Pol^3(m)$ $(m>4)$ with  nontrivial colorings is exactly
$$N_{nt}(m)=\begin{cases}
\sum_{0\leq k\leq {m\over 2}}([{k\over 2}]+1) & \text{ if  $m$ is even } \\
\sum_{1\leq
k\leq {m\over 2}}([{k\over 2}]+1)  & \text{ if  $m$ is odd.}\\
 \end{cases}$$
\end{cor}

 Next let us look at the case in which $\lambda$ is
trivial. By Lemma~\ref{alg-2} we divide our argument into two cases:
(I) $\Delta(\lambda)$ is odd; (II) $\Delta(\lambda)$ is even.

\vskip .2cm

{\bf Case (I): $\Delta(\lambda)$ is odd}.
\begin{lem}\label{alg-4}
Let $\lambda$ be trivial such that $\Delta(\lambda)$ is odd. Then
$$\bar{\mathcal{B}}(\lambda)=\begin{cases}
(0,2^{m-2}-1) & \text{if $\lambda\approx \lambda_{C_1}$}\\
(0,2^{m-3}-1) & \text{if $\lambda\approx\lambda_{C_3}$}\\
(2^{m-4}-1,0) & \text{if $\lambda\approx\lambda_{C_8}$}\\
(2^{m-3}-1,0) & \text{if $\lambda\approx\lambda_{C_9}$}\\
(2^{m-4}-1,0) & \text{if $\lambda\approx\lambda_{C_{10}}$}
\end{cases}$$
\end{lem}
\begin{proof}
If $\lambda\approx \lambda_{C_1}$, using Lemma~\ref{alg-2} and by
direct calculations, one has that
$s_1s_2=s_2s_3=\cdots=s_{m-1}s_m=s_ms_1$, so $A_0$ may be written as
follows:
\begin{equation*}
\begin{bmatrix}
0 &   x_2    &  x_3 & x_4 & \cdots  & x_{m-2} & x_{m-1} \\
x_2  & 0 & x_1  &  0 & \cdots   & 0 & 0  \\
  x_3   & x_1 & 0 & x_1 & \cdots & 0 & 0\\
   x_4   & 0 & x_1 & 0 & \cdots & 0 &0\\
         \hdotsfor{7}          \\
  x_{m-2}  & 0 &  0  & 0&  \cdots & 0  &x_1\\
   x_{m-1}  & 0 &  0  & 0& \cdots & x_1  &0
\end{bmatrix}
\end{equation*}

\noindent where $x_1=[s_3s_4]$ and $x_i=[fs_{i+1}], i=2, ..., m-1$.
We see easily that $\{x_1,..., x_{m-1}\}$ forms a basis of
$\mathcal{K}_\lambda/(\mathcal{K}_\lambda\cap\mathcal{H}_\lambda^2)$
so
$\dim\mathcal{K}_\lambda/(\mathcal{K}_\lambda\cap\mathcal{H}_\lambda^2)=m-1$.
Then one may conclude that
$\bar{\mathcal{B}}(\lambda)=(0,2^{m-2}-1)$.

\vskip .2cm In a similar way as above, if $\lambda\approx
\lambda_{C_3}$, one has that
$[s_2s_3]=[s_3s_4]=[s_4s_5]=\cdots=[s_{m-1}s_m]=0$, so $A_0$ may be
written as follows:
\begin{equation*}
\begin{bmatrix}
0 &   x_3    &  x_4 & x_5 & \cdots & x_{m-3} &  \sum_{j\text{ is odd}}x_j  &  x_1+\sum_{j\text{ is even}}x_j  \\
x_1  & 0 & 0  &  0 & \cdots   & 0 & 0 &0 \\
  x_2   & 0 & 0 & 0 & \cdots & 0 & 0 &0\\
   x_3   & 0 & 0 & 0 & \cdots & 0 &0&0\\
         \hdotsfor{8}          \\
  x_{m-3}  & 0 &  0  & 0&  \cdots & 0 & 0 &0\\
   \sum_{j\text{ is odd}}x_j  & 0 &  0  & 0& \cdots & 0&0  &0
\end{bmatrix}
\end{equation*}
 where $x_i=[fs_{i+1}], i=1, ..., m-3$. And $\{x_1, ...,
x_{m-3}\}$ forms a basis of
$\mathcal{K}_\lambda/(\mathcal{K}_\lambda\cap\mathcal{H}_\lambda^2)$
so
$\dim\mathcal{K}_\lambda/(\mathcal{K}_\lambda\cap\mathcal{H}_\lambda^2)=m-3$.
A direct observation shows that
$\bar{\mathcal{B}}(\lambda)=(0,2^{m-3}-1)$.

\vskip .2cm

If $\lambda\approx \lambda_{C_8}$ or $\lambda_{C_9}$, then
$[s_1s_2]=[s_2s_3]=\cdots=[s_{m-1}s_m]=[s_ms_1]=0$, so $A_0$ can be
reduced to a $1\times (m-3)$ matrix
$$([fs_2], [fs_4],  ..., [fs_{m-1}]).$$
Also, we easily see that $\{s_3^2, fs_2, fs_4, ..., fs_{m-1}\}$ can
be used as a basis of $\mathcal{K}_\lambda$ and $\{f^2, s_3^2\}$
forms a basis of $\mathcal{H}_\lambda^2$ (note that $\dim
\mathcal{H}_\lambda^2=m-1-\Delta(\lambda)$=2). However, when
$\lambda\approx \lambda_{C_8}$, by direct calculations one has that
$f^2=fs_2+\sum_{j>4\text{ is odd}}fs_j$,  so $f^2\in
\mathcal{K}_\lambda$ and $\mathcal{H}_\lambda^2\subset
\mathcal{K}_\lambda$. Thus, $\dim\mathcal{H}_\lambda^2\cap
\mathcal{K}_\lambda=2$,
$\dim\mathcal{K}_\lambda/(\mathcal{K}_\lambda\cap\mathcal{H}_\lambda^2)=m-4$
and $\{[fs_4],  ..., [fs_{m-1}]\}$ forms a basis of
$\mathcal{K}_\lambda/(\mathcal{K}_\lambda\cap\mathcal{H}_\lambda^2)$.
Moreover,  one has that
$\bar{\mathcal{B}}(\lambda_{C_8})=(2^{m-4}-1,0)$. When
$\lambda\approx \lambda_{C_9}$, it is not difficult to check that
$\dim\mathcal{H}_\lambda^2\cap \mathcal{K}_\lambda=1$ and $[fs_2],
[fs_4], ..., [fs_{m-1}]$ are linearly independent, so
$\dim\mathcal{K}_\lambda/(\mathcal{K}_\lambda\cap\mathcal{H}_\lambda^2)=m-3$.
Thus,  $\bar{\mathcal{B}}(\lambda_{C_9})=(2^{m-3}-1,0)$.

\vskip .2cm If $\lambda\approx \lambda_{C_{10}}$, then
$[s_1s_2]=[s_2s_3]=\cdots=[s_{m-1}s_m]=[s_ms_1]=[s_3^2]=0$ and
$fs_2=f^2$, so $A_0$ can be reduced to a $1\times (m-3)$ matrix $(0,
[fs_4], ..., [fs_{m-1}])$. It is easy to see that
$\dim\mathcal{H}_\lambda^2\cap \mathcal{K}_\lambda=1$ and
$\dim\mathcal{K}_\lambda/(\mathcal{K}_\lambda\cap\mathcal{H}_\lambda^2)=m-4$,
so $\bar{\mathcal{B}}(\lambda_{C_{10}})=(2^{m-4}-1,0)$.
\end{proof}

{\bf Case (II): $\Delta(\lambda)$ is even}.

\begin{lem}\label{alg-5}
Let $\lambda$ be trivial such that $\Delta(\lambda)$ is even. Then
$$\bar{\mathcal{B}}(\lambda)=\begin{cases}
(1, 2^{m-2}-2)  & \text{if $\lambda\approx \lambda_{C_2}$}\\
(2^{m-2}-1, 0) & \text{if $\lambda\approx\lambda_{C_4}$}\\
(2^{m-3}-1, 0) & \text{if $\lambda\approx\lambda_{C_5}$}\\
(2^{m-4}-1, 0) & \text{if $\lambda\approx\lambda_{C_6}$}\\
(2^{m-3}-1, 0) & \text{if $\lambda\approx\lambda_{C_7}$}
\end{cases}$$
\end{lem}
\begin{proof}
If $\lambda\approx \lambda_{C_2}$, then one can obtain by
Lemma~\ref{alg-2} that
$[s_1s_2]=[s_2s_3]=\cdots=[s_{m-1}s_m]=[s_ms_1]=[s_1^2]=0$ and
 so $A_0$ can be reduced to the following matrix

\begin{equation*}
\begin{bmatrix}
0 &   [fs_2]    &  [fs_4] &  \cdots   & [fs_{m-1}] \\
[fs_2]  & 0 &   0 & \cdots    & 0  \\
  [fs_3]    & 0 & 0 & \cdots  & 0\\
  [fs_4]   & 0 &  0 & \cdots  &0\\
         \hdotsfor{5}          \\
     [fs_{m-1}]  & 0 &  0 & \cdots   &0
\end{bmatrix}
\end{equation*}
 One may easily show that $\{[fs_2], ..., [fs_{m-1}]\}$ is
a basis of
$\mathcal{K}_\lambda/(\mathcal{K}_\lambda\cap\mathcal{H}_\lambda^2)$.
Then a direct observation can obtain that
$\bar{\mathcal{B}}(\lambda_{C_2})=(1, 2^{m-2}-2)$.

\vskip .2cm If $\lambda\approx \lambda_{C_4}$, then one has that
$[s_1s_2]=[s_2s_3]=\cdots=[s_{m-1}s_m]$, so $A_0$ can be reduced to
the following matrix
\begin{equation*}
\begin{bmatrix}
[fs_3] &  [fs_4]    &  [fs_5] &  \cdots &[fs_{m-3}] & [fs_{m-2}] & [fs_{m-1}] & [fs_{m}] \\
[s_1s_2]  & 0 & 0  &   \cdots  &0 & 0 & 0  &0\\
 0   & [s_1s_2] & 0 &  \cdots & 0&0 & 0&0\\
   [s_1s_2]   & 0 & [s_1s_2] &  \cdots &0& 0 &0&0\\
         \hdotsfor{8}          \\
   0  & 0 &  0  &  \cdots & [s_1s_2] & 0  & [s_1s_2]& 0\\
     0  & 0 &  0  &  \cdots&0 & [s_1s_2] & 0  & [s_1s_2]
\end{bmatrix}
\end{equation*}
 and $\{[s_1s_2],[fs_3], ...,  [fs_{m}]\}$ is a basis of
$\mathcal{K}_\lambda/(\mathcal{K}_\lambda\cap\mathcal{H}_\lambda^2)$
so
$\dim\mathcal{K}_\lambda/(\mathcal{K}_\lambda\cap\mathcal{H}_\lambda^2)=m-1$.
Furthermore, one knows that $\bar{\mathcal{B}}(\lambda_{C_4})=(
2^{m-2}-1, 0)$.

\vskip .2cm If $\lambda\approx \lambda_{C_5}$, then one has that
$[s_1^2]=[s_1s_2]=[s_2s_3]=\cdots=[s_{m-1}s_m]=0$, so $A_0$ can be
reduced to a $1\times (m-3)$ matrix $([fs_3], [fs_4],  ...,
[fs_{m-1}])$, and $\{[fs_3], [fs_4],  ..., [fs_{m-1}]\}$ is a basis
of
$\mathcal{K}_\lambda/(\mathcal{K}_\lambda\cap\mathcal{H}_\lambda^2)$
so
$\dim\mathcal{K}_\lambda/(\mathcal{K}_\lambda\cap\mathcal{H}_\lambda^2)=m-3$.
Note that in this case $\dim\mathcal{H}_\lambda^2\cap
\mathcal{K}_\lambda=1$.  Thus,
$\bar{\mathcal{B}}(\lambda_{C_5})=(2^{m-3}-1, 0)$.

\vskip .2cm

If $\lambda\approx \lambda_{C_6}$ or $\lambda_{C_7}$, similarly to
the case $\lambda\approx \lambda_{C_5}$,  then one has that
$[s_1^2]=[s_1s_2]=[s_2s_3]=\cdots=[s_{m-1}s_m]=0$, so $A_0$ can be
reduced to a $1\times (m-3)$ matrix $$([fs_3], [fs_4],  ...,
[fs_{m-1}]).$$ As in the proof of cases $\lambda\approx
\lambda_{C_8}$ or $\lambda_{C_9}$,  we  see that $\{s_2^2, fs_3,
fs_4, ..., fs_{m-1}\}$ can be used as a basis of
$\mathcal{K}_\lambda$ and $\{f^2, s_2^2\}$ forms a basis of
$\mathcal{H}_\lambda^2$. However, when $\lambda\approx
\lambda_{C_6}$, it is easy to check that
$\mathcal{H}_\lambda^2\subset \mathcal{K}_\lambda$, so
$\dim\mathcal{H}_\lambda^2\cap \mathcal{K}_\lambda=2$ and
$\dim\mathcal{K}_\lambda/(\mathcal{K}_\lambda\cap\mathcal{H}_\lambda^2)=m-4$.
Moreover, $\bar{\mathcal{B}}(\lambda_{C_6})=(2^{m-4}-1, 0)$ and
$b_{\Delta(\lambda)}=0$. When $\lambda\approx \lambda_{C_7}$, one
may check that $\dim\mathcal{K}_\lambda\cap\mathcal{H}_\lambda^2=1$
and then
$\dim\mathcal{K}_\lambda/(\mathcal{K}_\lambda\cap\mathcal{H}_\lambda^2)=m-3$.
Thus  $\bar{\mathcal{B}}(\lambda_{C_7})=(2^{m-3}-1, 0)$.
\end{proof}

\begin{rem} \label{5,6}
We see that for   $\lambda_{C_5}$ and  $\lambda_{C_7}$,
$\Delta(\lambda_{C_5})=\Delta(\lambda_{C_7})$ and
$\bar{\mathcal{B}}(\lambda_{C_5})=\bar{\mathcal{B}}(\lambda_{C_7})$.
However, we can still distinguish them by using the first
Stiefel-Whitney class. Let $w_1(\lambda)\in H^1(M(\lambda);\Z_2)$
denote the first Stiefel-Whitney class. It is well-known that
$w_1(\lambda)=0$ if and only if $M(\lambda)$ is orientable (see,
e.g. \cite{ms}). Then, by Corollary~\ref{ori-1} one knows that if
$\lambda\approx \lambda_{C_5}$, then $w_1(\lambda_{C_5})=0$; but  if
$\lambda\approx \lambda_{C_7}$, $w_1(\lambda_{C_7})\not=0$. This
also happens for $\lambda_{C_8}$ and $\lambda_{C_{10}}$, i.e., $\Delta(\lambda_{C_8})=\Delta(\lambda_{C_{10}})$ and
$\bar{\mathcal{B}}(\lambda_{C_8})=\bar{\mathcal{B}}(\lambda_{C_{10}})$. In this case we can
use the number $\dim\mathcal{K}_\lambda\cap\mathcal{H}_\lambda^2$ to
distinguish $\lambda_{C_8}$ and $\lambda_{C_{10}}$. Actually, by Lemma~\ref{alg-4}, if $\lambda\approx
\lambda_{C_8}$,
$\dim\mathcal{K}_\lambda\cap\mathcal{H}_\lambda^2=2$; but  if
$\lambda\approx \lambda_{C_{10}}$,
$\dim\mathcal{K}_\lambda\cap\mathcal{H}_\lambda^2=1$.
\end{rem}

\begin{thm}\label{c-r2}
Let $\lambda_1, \lambda_2$ be two trivial colorings on $\Pol^3(m)$
with $m>4$. Then $M(\lambda_1)$ and $M(\lambda_2)$ are homeomorphic
if and only if  their cohomologies $H^*(M(\lambda_1);\Z_2)$ and
$H^*(M(\lambda_2);\Z_2)$ are isomorphic as rings.
\end{thm}

\begin{proof}
This follows immediately from Lemmas~\ref{alg-2},
\ref{alg-4}-\ref{alg-5} and Remark~\ref{5,6}.
\end{proof}

As a consequence of Corollary~\ref{nu} and Theorem~\ref{c-r2}, one
has

\begin{cor}\label{number-2}
The number  of homeomorphism classes of small covers over $\Pol^3(m)
(m>4)$ with  trivial colorings is exactly
$$N_t(m)=\begin{cases}
4  & \text{ if  $m$ is odd}\\
6  &\text{ if  $m$ is even} \end{cases}$$
\end{cor}

\section{Proofs of Theorems~\ref{RIGIDITY}
and~\ref{NUMBER}}\label{proof}

Now let us finish the proofs of Theorems~\ref{RIGIDITY}
and~\ref{NUMBER}.

\vskip .2cm

\noindent
{\em Proof of Theorem~\ref{RIGIDITY}}. It suffices to show that if
 $H^*(M(\lambda_1);{\Bbb Z}_2)$ and
$H^*(M(\lambda_2);{\Bbb Z}_2)$ are isomorphic as rings, then
$M(\lambda_1)$ and $M(\lambda_2)$  are homeomorphic. By
Propositions~\ref{m>6}, \ref{c-r1} and~\ref{c-r2}, this is true when
$m>6$. It remains to consider the case $m\leq 6$. As stated in
Remark~\ref{3,4}, the cohomological rigidity holds when $m\leq 4$
(see also \cite{ly} and \cite{cmo}). Next, we only need pay our
attention on  the case $5\leq m\leq 6$. By Lemmas~\ref{alg-1},
\ref{alg-2}, \ref{alg-3}, \ref{alg-4}, \ref{alg-5} and
Remark~\ref{5,6},  we may list all possible $\lambda$ with mentioned
invariants in the case $5\leq m\leq 6$ whichever $\lambda$ is
trivial or nontrivial.

\vskip .2cm

(A) Case $m=5$:
\begin{center}
\begin{tabular}{|c|c|c|c|c|c|}
\hline &&&&&\\
 $\lambda$ & Trivialization & $\Delta(\lambda)$ &
$\bar{\mathcal{B}}(\lambda)$ & $(n_\lambda, m_\lambda)$ &
$w_1(\lambda)$\\
\hline
  $\lambda_{C_*}$   &   nontrivial  & 1 & $(1,0)$ &  $(1, 0)$ &\\
\hline
 $\lambda_{C_*}$     &  nontrivial  & 2 & $(1, 3)$  & $(2,0)$ &\\
\hline
  $\lambda_{C_*}$     &nontrivial     &  1 & $(1,0)$ & $(2, 2)$&\\
\hline
$\lambda_{C_3}$  &  trivial  & 3 & $(0, 3)$ & &\\
\hline
$\lambda_{C_5}$  & trivial & 2 & $(3, 0)$ &&0\\
\hline
$\lambda_{C_6}$& trivial & 2 & $(1, 0)$ &&\\
\hline
$\lambda_{C_7}$ & trivial &  2 & $(3,0)$ && nonzero\\
\hline
 \end{tabular}
\end{center}

\vskip .2cm

(B) Case $m=6$:
\begin{center}
\begin{tabular}{|c|c|c|c|c|c|}
\hline &&&&&\\ $\lambda$  & Trivialization &   $\Delta(\lambda)$ &
$\bar{\mathcal{B}}(\lambda)$ & $(n_\lambda, m_\lambda)$ &
$\dim\mathcal{K}_\lambda\cap\mathcal{H}_\lambda^2$\\
\hline
 $\lambda_{C_*}$  &      nontrivial  & 1 & $(1,0)$ &  $(0, 0)$ &\\
     \hline
  $\lambda_{C_*}$  &    nontrivial     &  1 & $(1,0)$ & $(1,0)$&\\
\hline
$\lambda_{C_*}$  &       nontrivial  & 2 & $(1, 3)$  & $(2,0)$ &\\
     \hline
 $\lambda_{C_*}$  &     nontrivial     &  3 & $(0,3)$ & $(3,0)$&\\
\hline
  $\lambda_{C_*}$  &    nontrivial     &  1 & $(1,0)$ & $(2, 2)$&\\
    \hline
$\lambda_{C_*}$  &      nontrivial     &  2 & $(1,1)$ & $(3, 2)$&\\
\hline
$\lambda_{C_1}$  &  trivial  & 5 & $(0, 15)$ & &\\
\hline
$\lambda_{C_2}$  &    trivial & 4 & $(1,14)$ &&\\
\hline
$\lambda_{C_4}$  &   trivial & 4 & $(15, 0)$ &&\\
\hline
 $\lambda_{C_8}$  &   trivial & 3 & $(3,0)$ & &2\\
  \hline
 $\lambda_{C_9}$  &   trivial & 3 & $(7,0)$ &&\\
  \hline
 $\lambda_{C_{10}}$  &   trivial & 3 & $(3,0)$ & &1\\
\hline
 \end{tabular}
\end{center}
We clearly see from two tables above that by using invariants
$\Delta(\lambda)$, $\bar{\mathcal{B}}(\lambda)$, $(n_\lambda,
m_\lambda)$, $w_1(\lambda)$ and
$\dim\mathcal{K}_\lambda\cap\mathcal{H}_\lambda^2$, we can
distinguish all $M(\lambda)$ up to homeomorphism when $m=5,6$. This
completes the proof. $\Box$

\vskip .2cm

Furthermore, Theorem~\ref{NUMBER} follows immediately  from
Theorem~\ref{RIGIDITY}, Corollaries~\ref{number-1}, \ref{number-2}
and Remark~\ref{3,4}.

\vskip .2cm

Finally, let us return to the invariants $\Delta(\lambda)$ and
$\mathcal{B}(\lambda)$ again. We see that generally these invariants
can always be defined for any small cover over a simple convex
polytop $P^n$. We would like to pose the following problems:
\begin{enumerate}
\item[$\bullet$]Under what condition can $\Delta(\lambda)$ and
$\mathcal{B}(\lambda)$ become the combinatorial invariants?
\item[$\bullet$]If $\Delta(\lambda)$ and
$\mathcal{B}(\lambda)$ are  the combinatorial invariants, then  how
can one calculate them in terms of polytopes $P^n$?
\end{enumerate}

\begin{acknow} The authors would like to express their thanks to the referee for valuable comments and suggestions.
\end{acknow}

\end{document}